\documentclass[11pt,reqno,oneside]{amsart}
\usepackage{amsmath,amsthm,amssymb,amsfonts,amsxtra,comment,graphicx,psfrag,mathrsfs}
\usepackage{mathtools,color,hhline,pdfsync}
\usepackage{hyperref,enumerate,enumitem,bm}
\usepackage[margin=1.0in]{geometry}
\usepackage{fancyhdr}
\allowdisplaybreaks

\allowdisplaybreaks

\def\R{{\mathbb R}}

\def\d{{\mathrm d}}

\theoremstyle{plain}
\newtheorem{theorem}{Theorem}[section]
\newtheorem{lemma}[theorem]{Lemma}
\newtheorem{corollary}[theorem]{Corollary}

\numberwithin{equation}{section}
\numberwithin{table}{section}
\numberwithin{figure}{section}
\usepackage{titlesec} 
\titleformat{\section}{\vskip10pt\normalsize\bfseries}{\thesection.}{0.5em}{\centering}
\titleformat{\subsection}{\vskip10pt\normalsize\bfseries}{\thesubsection.}{0.5em}{}
\def\R{{\mathbb R}}

\def\d{{\mathrm d}}

\newcommand{\bR}{\mathbf{R}}
\newcommand{\bP}{\mathbf{P}}
\newcommand{\nn}{\nonumber}
\def\d{{\rm d}}

\begin{document}

\title{Analysis of fully discrete FEM for 
miscible displacement in porous media with
Bear--Scheidegger diffusion tensor} 

\author[Wentao Cai]{Wentao Cai}
\address{Department of Mathematics, School of Sciences, Hangzhou Dianzi University, Hangzhou 310018, China} 
\email {femwentao@hdu.edu.cn} 

\author[Buyang Li]{Buyang Li}
\address{Department of Applied Mathematics, 
The Hong Kong Polytechnic University, Kowloon, Hong Kong} 
\email {bygli@polyu.edu.hk} 

\author[Yanping Lin]{Yanping Lin}
\address{Department of Applied Mathematics, The Hong Kong Polytechnic University, Kowloon, Hong Kong. E-mail: yanping.lin@ polyu.edu.hk} 
\email {yanping.lin@polyu.edu.hk} 

\author[Weiwei Sun]{Weiwei Sun}
\address{Department of Mathematics, City University of Hong Kong,  
Kowloon, Hong Kong.}
\email {maweiw@cityu.edu.hk}

\setlength\abovedisplayskip{3pt}
\setlength\belowdisplayskip{3pt}

\date{} 

\begin{abstract}
{
Fully discrete Galerkin finite element methods are studied 
for the equations of miscible displacement in porous media with the 
commonly-used Bear--Scheidegger diffusion-dispersion tensor: 
$$
D({\bf u}) = \gamma d_m I + |{\bf u}|\bigg( \alpha_T I +
(\alpha_L - \alpha_T) \frac{{\bf u} \otimes {\bf u}}{|{\bf u}|^2}\bigg) \, .
$$
Previous works on optimal-order $L^\infty(0,T;L^2)$-norm error estimate 
required the regularity assumption  
$\nabla_x\partial_tD({\bf u}(x,t)) \in L^\infty(0,T;L^\infty(\Omega))$,   
while the Bear--Scheidegger diffusion-dispersion tensor is only Lipschitz 
continuous even for a smooth velocity field ${\bf u}$.
In terms of  the maximal $L^p$-regularity of fully discrete finite element 
solutions of parabolic equations, optimal error estimate 
in $L^p(0,T;L^q)$-norm and 
almost optimal  error estimate in $L^\infty(0,T;L^q)$-norm are established 
under the assumption of $D({\bf u})$ being Lipschitz continuous 
with respect to ${\bf u}$.\\

\noindent{\bf Keywords.} 
miscible displacement in porous media,
Bear--Scheidegger diffusion-dispersion tensor,
finite element method,
maximal $L^p$-regularity, error estimate
}
\end{abstract}

\maketitle

\section{Introduction}\label{Sec:intr}
The incompressible flow of binary miscible fluid in porous media 
is governed by the miscible displacement equations
\begin{align}
& \gamma \frac{\partial{c}}{\partial{t}}-\nabla\cdot(D({\bf{u}})\nabla c)
+{\bf{u}}\cdot\nabla c =\hat{c}q_{I}-cq_{I},
\label{e1} \\
& \nabla\cdot{\bf u} =q_I-q_P, \qquad  {\bf u} =-\frac{k(x)}{\mu(c)}\nabla p , 
\label{e2}
\end{align}
where ${\bf u}$ and $p$ are the velocity and pressure of the fluids mixture,  
respectively, and $c$ is the concentration of one fluid. In this model,  
$k(x)$ is the permeability of the porous medium, $\mu(c)$ the  
concentration-dependent viscosity, $\gamma$ the porosity of the medium,  
$q_I\ge 0$ and $q_P\ge 0$ the given injection and production sources,  
respectively, 
and $\hat{c}$ the concentration in the injection source. A popular  
diffusion-dispersion tensor $D({\bf {u}})=[D_{ij}({\bf{u}})]_{d \times d}$  
used in reservoir simulations and underground oil exploration 
is the Bear--Scheidegger model (cf. \cite{BB2,Sche})
\begin{equation}\label{expr-Du}
D({\bf u}) = \gamma d_m I + |{\bf u}| \left ( \alpha_T I +
(\alpha_L - \alpha_T) \frac{{\bf u} \otimes {\bf u}}{|{\bf u}|^2}
\right ) ,
\end{equation}
where 
$d_m>0$ denotes the molecular diffusion, and 
$\alpha_L$ and $\alpha_T$ the constant longitudinal and transversal 
dispersivities of the isotropic porous medium, respectively.
We consider \eqref{e1}-\eqref{e2} in a bounded smooth domain 
$\Omega\subset{\mathbb R}^d$, with $d\in\{2,3\}$, up to time $T$,  
subject to the no-flux boundary conditions
\begin{align}
&{\bf{u}}\cdot{\bf {n}}=0
\quad\mbox{and}\quad
D({\bf{u}}) \nabla c\cdot {\bf {n}}=0
\quad \mbox{on}\,\,\, \partial\Omega\times(0,T] , \label{b1}
\end{align}
with the given initial condition
\begin{align}
&c(x, 0)=c_0(x)\quad\mbox{for}\,\,\, x\in \Omega   .
\label{b2}
\end{align}

Numerical methods and analysis for the miscible displacement system  
(\ref{e1})-(\ref{b2}) have been investigated extensively in the last  
several decades, and numerical simulations have been done for various  
engineering applications, $e.g.$, \cite{DEW, EW, Feng,
Wang, WLELQ, WSS}.  A traditional approach to establish the optimal  
$L^\infty(0,T;L^2)$-norm error estimate is based on an elliptic Ritz  
projection
$\bR_h(t):H^1(\Omega)\rightarrow S_h^r$ onto the finite element space, defined by (see \cite{Whe})
\begin{align}
\Big(D({\bf u}(\cdot,t))\nabla (\phi-\bR_h\phi), \, \nabla \varphi_h \Big)
= 0 ,\quad
\mbox{for~all}~~\phi\in H^1(\Omega)~~\mbox{and}~~\varphi_h\in S_h^r  \, .
\end{align}
Most previous works on optimal $L^\infty(0,T;L^2)$ error estimates of
Galerkin type FEMs for \eqref{e1}-\eqref{b2} follow this way, 
which requires the following estimate of the Ritz projection:
\begin{align}
\label{p-est}
\|\partial_t(c - \bR_hc)\|_{L^2(0,T;L^2)}\leq Ch^{r+1} \, .
\end{align}
The estimate above was established by Wheeler \cite{Whe} 
under the regularity assumption
\begin{equation}
\|\nabla_x\partial_t D({\bf u}(x,t))\|_{L^\infty(0,T;L^\infty)}\leq C
\label{Dxt}
\end{equation}
for a general nonlinear parabolic equation. However, less attention was  
paid to the regularity of the Bear--Scheidegger diffusion--dispersion  
tensor. It was shown in \cite{SW} that
$D({\bf u})$ is Lipschitz continuous in ${\bf u}$. In a more recent work
\cite{LS2}, a counter example was presented to show that
even for a smooth velocity field it may hold
$$
\nabla_x\partial_tD({\bf u}(x,t))
\notin L^p(\Omega_T)
\quad\mbox{for any}~~ p\ge 1 .
$$
Clearly, the Bear--Scheidegger dispersion model may not satisfy  
the regularity condition \eqref{Dxt} and therefore, optimal  
$L^\infty(0,T;L^2)$ error estimates of fully discrete Galerkin-Galerkin FEMs,  
Galerkin-mixed FEMs and many other numerical methods for this model have  
not been well investigated in this case.

In this article, we study the commonly-used Bear--Scheidegger 
diffusion-dispersion  
model by a linearized fully discrete Galerkin FEM and establish an  
optimal $L^p(0,T;L^q)$ error estimate, together with an almost optimal  
$L^\infty(0,T;L^q)$ error estimate. The key to our analysis is the discrete  
maximal $L^p$-regularity ($L^p$-stability) of fully discrete finite element  
solutions of the parabolic equations 
\begin{align}\label{leq}
\left\{
\begin{array}{ll}
\partial_t\phi - \nabla \cdot (a\nabla \phi ) +\phi= f -\nabla\cdot{\bf g}
&\mbox{in}~\Omega,
\\[5pt]
\displaystyle
a\nabla\phi\cdot{\bf n}
={\bf g} \cdot{\bf n}
&\mbox{on} 
~~\partial\Omega ,\\[5pt]
\phi(x,0)=\phi_0(x) &\mbox{for}~
x\in \Omega \, . 
\end{array}
\right.
\end{align}
In the last several decades,  great efforts have been devoted to the 
maximal $L^p$-stability estimates, $e.g.$, see 
\cite{Cro,DLSW,Gei1,Gei2,Han,Ley,LV1,Lin2, LTW, NW,Pal,STW1,STW2} and  
references therein. A straightforward application of the maximal  
$L^p$-stability estimates is the error estimates 
\begin{align}
&\| \bP_h \phi  - \phi_h \|_{L^p(0,T; L^q)}
\le C (\| \bP_h \phi_0(x) - \phi_h(0)\|_{L^q}
+ \|\phi - \bR_h\phi \|_{L^p(0,T; L^q)}) ,\label{error}
  \\
&\| \bP_h \phi - \phi_h \|_{L^\infty(0,T; L^q)}
\le C\|\bP_h \phi_0(x) - \phi_h(0)\|_{L^q}
+ C\ln(2+1/h)\| \bP_h \phi - \phi \|_{L^\infty(0,T; L^q)}, 
\label{error2}
\end{align}
with $p,q \in(1,\infty)$, where $\phi_h$ is the finite element solution of \eqref{leq}, $\bP_h$ is the $L^2$-projection operator 
onto the finite element space $S_h^r$, and $\bR_h$ the Ritz projection  
operator associated with the elliptic operator  
$\mathcal{L}=- \nabla \cdot (a\nabla) +1 $. 
Early works on  
such $L^p(0,T;L^q)$ and $L^\infty(0,T;L^q)$ stability estimates  
were done mainly for spatially semi-discrete finite element solutions  
of linear parabolic equations with sufficiently smooth time-independent  
coefficients, $e.g.$, $a_{ij} = a_{ij}(x) \in C^{2+\alpha}(\overline\Omega)$. 
The extension to time-independent Lipschitz continuous coefficients  
$a_{ij}=a_{ij}(x) \in W^{1,\infty}(\Omega)$ was presented in \cite{Li}.   
Further extensions to fully discrete finite element solutions were  
done in \cite{KS,KLL,LV1} for linear autonomous parabolic equations and  
in \cite{LS3} for linear nonautonomous parabolic equations  
(with coefficients $a_{ij}=a_{ij}(x,t)$). 
The former relies on  
the semigroup approach which is applicable only for a problem with  
time-independent coefficients, and the latter uses a perturbation  
technique together with a duality argument.

The $L^p(0,T; L^q)$ approach 
has apparent advantages over the traditional  
$L^\infty(0,T;L^2)$ estimate in dealing with nonlinear parabolic equations.  
Recently, analysis on semi-discrete nonlinear parabolic equations was 
presented by several authors, see \cite{Gei2,LS2} for semi-discrete finite 
element methods and \cite{AkrivisLi,ALL,KuLL} for time discrete systems. 
However, 
no analysis has been done for fully discrete Galerkin FEMs  
for nonlinear physical equations. 
The $L^p(0,T; L^q)$ analysis of  
a fully discrete FEM for nonlinear parabolic equations is much different from  
the analysis of time-discrete systems. 
In this paper, we apply the $L^p(0,T; L^q)$ approach to  
commonly-used linearized fully discrete Galerkin finite element methods for 
the nonlinear miscible displacement problem (\ref{e1})-(\ref{b2}) 
with the Bear--Scheidegger diffusion-dispersion tensor to establish 
optimal $L^p(0,T; L^q)$ and almost optimal $L^\infty(0,T; L^q)$ 
error estimates.  
More important is that our analysis illustrates a fundamental tool in establishing optimal 
error estimates of commonly-used fully discrete Galerkin FEMs for 
nonlinear physical equations with more general diffusion coefficients.

%
%
%
%
\section{Main results}
\setcounter{equation}{0}
For $q\in[1,\infty]$ and any integer $k\ge 0$, we denote by  
$W^{k,q}=W^{k,q}(\Omega)$ the usual Sobolev spaces of functions  
defined on $\Omega$, with the abbreviations $L^q=W^{0,q}$ and  
$H^k=W^{k,2}$; see \cite{Adams}. The dual  
space of $W^{k, q}$ is denoted by $\widetilde{W}^{-k, q'}$,  
with the notation $q'=q/(q-1)$ and the abbreviation  
$\widetilde{H}^{-k}=\widetilde{W}^{-k, 2}$. For any integer $k\ge 0$  
and $\alpha\in(0,1)$, we denote by  $C^{k,\alpha}$ 
the space of functions whose partial derivatives up to $k^{\rm th}$-order  
are H\"older continuous with the exponent $\alpha$.

Let $0=t_0<t_1<\cdot\cdot\cdot<t_N=T$ be a uniform partition of the  
interval $[0, T]$ for some integer $N$, with the step size  
$t_n-t_{n-1}=\tau=T/N$. For any sequence of functions $\{f^n\}^N_{n=0}$,  
we define
\begin{align*}
&D_{\tau}f^n:=\frac{f^n-f^{n-1}}{\tau}, \\[5pt]
&\| f^m \|_{L^p (X)}: =
\left\{ 
\begin{aligned}
&\mbox{$\left ( \sum_{n=1}^{m} \tau \|  f^n \|_X^p \right )^{\frac1p} $} ,   
&& p\in[1,\infty), \\ 
& \max_{1\le n\le m}  \|f^n\|_X   ,  && p=\infty , 
\end{aligned} 
\right.
\end{align*}
for certain Sobolev space $X$. 
The norm $\| f^m \|_{L^p (X)}$ is simply the $L^p(0,m\tau;X)$ norm of the piecewise constant function which takes the value $f_n$ on each interval $(t_{n-1},t_n]$. 

Let $\Omega\subset{\mathbb R}^d$, with $d\in\{2,3\}$, be a bounded domain with smooth boundary $\partial\Omega$, and let $\mathcal{T}_h$ be a shape-regular and quasi-uniform triangulation of $\Omega$ into triangles or tetrahedra which fit the boundary $\partial\Omega$ exactly, with possibly curved triangles or tetrahedra near on the boundary. 
We denote by $h$ the mesh size of triangulation, and define the following finite element spaces: 
\begin{align*}
&\text {$S^r_h=\{\phi_h\in H^1(\Omega)$: $\phi_h$ is a polynomial  
of degree $r$ on each triangle (or tetrahedra)\}}, 
\\
&\mathring S^{2}_{h}=\{\phi_h\in S^{2}_{h}:  
\mbox{$\int_\Omega\phi_h \d x=0$}\}.
\end{align*}
We consider a linearized and stabilized fully-discrete FEM 
for (\ref{e1})-(\ref{b2}), which seeks $P^{n-1}_h\in \mathring S^{2}_{h}$  
and $\mathcal{C}^{n}_h\in S^1_h$ such that
\begin{align}
& \left({\frac{k(x)}{\mu(\mathcal{C}_h^{n-1})}}\nabla P_h^{n-1},  
\nabla v_h\right) =(q_I^{n-1}-q_P^{n-1}, v_h),
\quad\forall\, v_h\in  \mathring S^{2}_h , \quad n=1,\dots,N+1,
\label{fdis2}\\[10pt]
& (\gamma D_{\tau}\mathcal{C}_h^{n}, w_h)
+(D({\bf {U}}_h^{n-1})\nabla \mathcal{C}_h^{n}, \nabla w_h)
+\left(\frac{1}{2}(q_I^n+q_P^n)\mathcal{C}^n_h, w_h\right)
\label{fdis1} \\
& +\frac{1}{2}({\bf{U}}^{n-1}_h\cdot\nabla \mathcal{C}^n_h, w_h)
-\frac{1}{2}({\bf{U}}^{n-1}_h \mathcal{C}^n_h, \nabla w_h)=(\hat{c}q_I^n, w_h),
\quad\forall\, w_h\in S^1_h,
\quad n=1,\dots,N,  \nonumber
\end{align}
where
\begin{align}
{\bf{U}}^{n-1}_h=-\frac{k(x)}{\mu(\mathcal{C}^{n-1}_h)}\nabla P^{n-1}_h   
,\label{TDIS4}
\end{align}
and $\mathcal{C}^0_h=\Pi_hc(\cdot,0)$, with $\Pi_h$ being the 
Lagrange interpolation operator onto $S_h^1$.

We assume that $q_I$, $q_P$, $\hat{c} \in C([0,T];L^\infty(\Omega))$,  
$k\in W^{2,\infty}(\Omega)$, $\mu\in W^{2,\infty}(\mathbb{R})$,  
$k_0\leq k(x)\leq k_1$, $\mu_0\leq \mu(c)\leq\mu_1$, and the system  
\eqref{e1}-\eqref{b2} 
has a unique solution satisfying
\begin{align}\label{regu}
 \|c\|_{C([0,T];W^{2,q})}+\|\partial_t c\|_{C([0,T];W^{1, q})} 
+\|\partial_{tt}c\|_{C([0,T];\widetilde{W}^{-1, q})} 
+\| p\|_{ C([0,T];W^{3,q})} \le K .
\end{align}
This only guarantees the Lipschitz continuity $D({\bf u})\in L^\infty(0,T;W^{1,\infty})\cap W^{1,\infty}(0,T;L^\infty)$, instead of \eqref{Dxt}, for the Bear--Scheidegger diffusion-dispersion tensor \eqref{expr-Du}.  
Our main result is presented in the following theorem, with the notations
\begin{align}
c^n=c(\cdot,t_n),
\quad
p^n=p(\cdot,t_n),
\quad\mbox{and}\quad
{\bf u}^n={\bf u}(\cdot,t_n) \, . 
\nn 
\end{align}

\begin{theorem} \label{MainTHM}
Suppose that the system \eqref{e1}-\eqref{b2} has a unique solution  
$(c, {\bf u}, p)$ satisfying \eqref{regu} for some $q\in(d,\infty)$.  
Then the finite element system \eqref{fdis2}-\eqref{TDIS4} 
admits a unique solution $(P^{n}_h,  \mathcal{C}^{n}_h)$, $n=1,\dots, N$, 
satisfying 
\begin{align}
\|P^{n}_h-p^{n}\|_{L^p(W^{1,q})}+\|{\bf{U}}^{n}_h-{\bf{u}}^{n}
\|_{L^p(L^q)} +\|\mathcal{C}_h^{n}-c^{n}\|_{L^p(L^q)}
\le C_{p,q}(\tau+h^2) ,  
\label{pqerror} 
\end{align} 
for any $p\in(1,\infty)$,  
where $C_{p,q}$ is a constant, independent of $n$, $\tau$
and $h$ and dependent upon $p,q$. 
\end{theorem}

\begin{corollary}\label{Corollary} 
Under the assumptions of Theorem \ref{MainTHM}, it holds that  
\begin{align} 
&\!\!\!\!\!\|P^{n}_h-p^{n}\|_{L^\infty(W^{1,q})}+\|{\bf{U}}^{n}_h 
-{\bf{u}}^{n} \|_{L^\infty(L^{q})} 
+\|\mathcal{C}_h^{n}-c^{n}\|_{L^\infty(L^q)}
\le C_{\epsilon}(\tau^{1-\epsilon}+h^{2-\epsilon}) , 
\label{Linftyerror}
\end{align}
for an arbitrary small $\epsilon>0$. 
\end{corollary}

The rest of this paper is devoted to the proofs of Theorem \ref{MainTHM} and Corollary \ref{Corollary}.  
The main difficulty iis to prove an upper bound for $\|\mathcal{C}_h^n\|_{W^{1,\infty}}$ in order to control the nonlinear terms involved in the analysis. To this end, we adopt the the error splitting approach developed in \cite{LiSun2013,LS13-SINUM-2} and the discrete maximal $L^p$-regularity of parabolic equations developed in \cite{KS,KLL,Li,LS2,LS3}. 
By this approach, we first prove in Section \ref{sec:lptime} that the semi-discretization in time has sufficient regularity uniformly with respect to the time-step size, $i.e.$, 
\begin{align*}
&\|D_{\tau}\mathcal{C}^{N}\|_{L^p(L^q)}+\|\mathcal{C}^{N}\|_{L^p(W^{2, q})}  
\le C_{p,q} ,
\end{align*}
where $C_{p,q}$ is a constant independent of the time-step size $\tau$. The estimate above implies an upper bound for $\|\mathcal{C}^n\|_{W^{1,\infty}}$ through the following discrete inhomogeneous Sobolev embedding: 
\begin{align*}
\|\mathcal{C}^{N} \|_{L^\infty(W^{1,\infty})}
\le 
C(\| D_{\tau}\mathcal{C}^{N}\|_{L^p(L^q)}+\|\mathcal{C}^{N}\|_{L^p(W^{2, q})} )
\le C_{p,q} ,
\end{align*}
which holds for sufficiently large $p$ and $q$ such that $\frac{2}{p}+\frac{d}{q}<1$. 

By using the regularity estimate above, in Section \ref{sec:lpfully}, we further derive error estimate for the fully discrete solution in the 
$L^p(\widetilde{W}^{-1, q})$ and $L^p(W^{1, q})$ norm, $i.e.$, 
\begin{align*}
\|D_{\tau}(\mathcal{C}^{n}_h-\Pi_h\mathcal{C}^{n})\|_ 
{L^p(\widetilde{W}^{-1, q})}
+\|\mathcal{C}^{n}_h-\Pi_h\mathcal{C}^{n}\|_{L^p(W^{1, q})}
\le Ch ,
\end{align*}
which yields an error estimate in $L^\infty(L^{\infty})$ through the discrete inhomogeneous Sobolev embedding 
\begin{align*}
\|\mathcal{C}^{n}_h-\Pi_h\mathcal{C}^{n}\|_{L^\infty(L^{\infty})}
\le 
C(\|D_{\tau}(\mathcal{C}^{n}_h-\Pi_h\mathcal{C}^{n})\|_ 
{L^p(\widetilde{W}^{-1, q})}
+\|\mathcal{C}^{n}_h-\Pi_h\mathcal{C}^{n}\|_{L^p(W^{1, q})})
\le Ch 
\end{align*}
for sufficiently large $p$ and $q$ such that $\frac{2}{p}+\frac{d}{q}<1$. 
By using the inverse inequality of the finite element space, we further obtain
\begin{align*}
\|\mathcal{C}^{n}_h-\Pi_h\mathcal{C}^{n}\|_{L^\infty(W^{1,\infty})}
\le 
Ch^{-1}\|\mathcal{C}^{n}_h-\Pi_h\mathcal{C}^{n} \|_{L^\infty(L^{\infty})}
\le C,
\end{align*}
which implies upper bound for $\|\mathcal{C}_h^n\|_{W^{1,\infty}}$.

Throughout we denote $C_{p_1,\dots,p_k}$ a generic positive constant  
which may be different at different  occurrence, independent of $n$, $\tau$  
and $h$, while possibly depend upon $K$, $T$, $\Omega$ and  the 
parameters $p_1,\dots,p_k$ in the subscript.


\section{Preliminaries} 
\label{Sec:Pre}
In this section we introduce some notations and lemmas to be used in our proof of Theorem \ref{MainTHM}. 
The basic ideas for proving these lemmas are described, and the detailed proof can be found in Appendix. 

We define a Ritz operator ${\bf {R}}_h(t):$ $H^1\rightarrow S^1_h$ and  
an $L^2$-projection operator ${\bf P}_h^r:L^2\rightarrow S_h^r$   
by 
\begin{align*}
(D({\bf u}(\cdot,t))\nabla(\phi-{\bf{R}}_h\phi), \nabla \varphi_h)  
+(\phi-{\bf{R}}_h\phi, \varphi_h)=0,
\quad\forall\,\phi\in H^1,\,\,\,\forall\, \varphi_h \in S^1_h ,
\end{align*}
and
\begin{align*}
(\phi-{\bf P}_h^r\phi, \varphi_h)=0,  \quad \forall\, \phi\in L^2 ,\,\,  
\forall\, \varphi_h\in S^r_h ,
\end{align*}
respectively, with the abbreviations ${\bf P}_h: ={\bf P}_h^1$ and  
$\overline{\bf P}_h:={\bf P}_h^2$,
which satisfy the following estimates:
\begin{align}
&\|\varphi-{\bf P}_h^r\varphi\|_{W^{\ell_0, q}} 
\leq Ch^{m-\ell_0}\|\varphi\|_{W^{m,q}},
&&\forall\varphi    \in W^{m,q},
\label{ellep}\\
&\|\varphi-{\bf{R}}_h\varphi\|_{L^s}+h\|\varphi
-{\bf{R}}_h\varphi\|_{W^{1,s}}\leq Ch^{l}\|\varphi\|_{W^{l,s}},
&&\forall\varphi    \in W^{l,s},
\label{ellep-2} \\
&\|\varphi-{\bf{R}}_h\varphi\|_{L^s}
\le Ch \|\varphi-{\bf{R}}_h\varphi\|_{W^{1,s}} ,
&&\forall\varphi    \in W^{l,s},
\label{ellep-33}
\end{align}
for $\ell_0=0, 1$, $\ell_0\leq m\leq r+1$, $1\le l \le r+1$, $1 
\leq q\leq\infty$ and $1<s<\infty$. 
Similarly, the Lagrangian  
interpolation operator $\Pi_h:C(\overline\Omega)\rightarrow S_h^1$ satisfies
\begin{align}\label{Lag-inter}
\|\Pi_h\varphi-\varphi\|_{L^q} 
+h\|\nabla(\Pi_h\varphi-\varphi)\|_{L^q}\leq C h^2\|\varphi\|_{W^{2,q}} , 
\quad  \forall\, \varphi\in W^{2,q},
\quad\forall\, q\in[2,\infty).
\end{align}

For the system \eqref{leq},
we define a corresponding time-discrete (spatially continuous) system
\begin{align}\label{leq-td}
\left\{
\begin{array}{ll}
D_{\tau} \Phi^n - \nabla \cdot (a(\cdot,t_n)\nabla \Phi^n ) + \Phi^n 
= f^n -\nabla\cdot{\bf g}^n
&\mbox{in}~\Omega,
\\[5pt]
\displaystyle
a(\cdot,t_n)\nabla\Phi^n\cdot{\bf n}
={\bf g}^n \cdot{\bf n}
&\mbox{on}
~~\partial\Omega ,\\[5pt]
\Phi^0=\phi_0(x) &\mbox{for}~
x\in \Omega ,
\end{array}
\right.
\quad n=1,\dots,N,
\end{align}
and a fully-discrete finite element system of  
$\Phi^n_{h}\in S_h^r$, $n=1,2,\dots$, 
\begin{align}\label{leq-fe}
(D_{\tau}\Phi^n_{h},v_h) + (a(\cdot,t_n)\nabla \Phi^n_h,\nabla v_h) 
+(\Phi^n_h,v_h)
= (f^n, v_h)+({\bf g}^n,\nabla v_h), \quad\forall\, v_h \in S_h^r  ,
\end{align}
where $f^n=f(\cdot,t_n)$ and ${\bf g}^n={\bf g}(\cdot,t_n)$. 
Some existing estimates for the solutions of \eqref{leq-td} and  
\eqref{leq-fe} are given in the following two lemmas.

\begin{lemma}
\label{lemma1}
{\it
If the coefficient matrix $a(x,t)=(a_{ij}(x,t))_{d\times d}$  
in \eqref{leq-td} and \eqref{leq-fe} satisfies
\begin{align}
&\lambda^{-1}\sum_{i=1}^d|\xi_i|^2\le
\sum_{i,j=1}^d a_{ij}(x,t)\xi_i\xi_j\le
\lambda \sum_{i=1}^d|\xi_i|^2 ,
\quad\forall\,\xi_i\in{\mathbb R},\,\,\, 
\forall\, (x,t)\in\Omega\times[0,T],   \label{ellipticity}\\
&a_{ij} \in L^\infty(0,T;W^{1,\infty}(\Omega))\quad\mbox{and}\quad
\partial_ta_{ij}  \in L^\infty(0,T;L^\infty(\Omega)) ,    
\label{W1infty-condition}
\end{align}
then the time-discrete solutions defined by \eqref{leq-td} satisfy
\begin{align}
&\|D_\tau \Phi^n \|_{L^p(\widetilde W^{-1,q})} + \| \Phi^n \|_{L^p(W^{1,q})}
\leq C (\|f^n\|_{L^p(L^q)}+\|{\bf g}^n\|_{L^p(L^q)}) ,&&\forall\,  
p,q\in(1,\infty) ,
\label{lemma-2}\\
&\|D_\tau \Phi^n\|_{L^p(L^q)} +\|\Phi^n \|_{L^p(W^{2,q})}\leq
 C \|f^n\|_{L^p(L^q)}, 
\qquad\mbox{if}\,\,\,{\bf g}={\bf 0},\,\, &&\forall\, p,q\in(1,\infty) \, . 
\label{lemma-1}
\end{align}
}
\end{lemma}

The proof of \eqref{lemma-1} was given in \cite{APW} (also see  
\cite[Theorem 3.1]{KLL}) and the proof for \eqref{lemma-2} can  
be found in \cite{LS3}. The following lemma is a consequence of \cite[(1.18)  
and (2.3)-(2.4)]{LS3}.

\begin{lemma}
\label{lemma2}
{\it 
Let $\phi^n=\phi(\cdot,t_n)$, $\Phi^n$ and $\Phi_h^n$ denote the solutions of \eqref{leq}, \eqref{leq-td} and \eqref{leq-fe}, respectively.  
Under the assumption of Lemma \ref{lemma1},
there exist positive constants $\tau_2$ and $h_2$ such that the  
following estimates hold 
for $\tau\le \tau_2$, $h \le h_2$ and $p,q \in (1, \infty):$
\begin{align}
&\!\!\! \|D_\tau \Phi_h^n\|_{L^p(\widetilde W^{-1,q})}+\|\Phi_h^n\|_{L^p(W^{1,q})}\leq
C \big(\|f^n\|_{L^p(L^q)} + \|{\bf g}^n\|_{L^p(L^q)}\big),
\label{lemma-3}
\\[5pt]
&\!\! \|{\bf P}_h \phi^n - \Phi_h^n \|_{L^p(L^q)}  \nonumber \\
&\le
C(\|\phi^n - {\bf R}_h \phi^n \|_{L^p(L^q)} 
+ \|{\bf P}_h \phi_0(x)-\Phi^0_h\|_{L^q}
+ \tau\|\partial_{tt}\phi \|_{L^p(0,T;\widetilde W^{-1,q})} ),
\label{lemma-4}
\\[5pt]
&\!\!\! \|D_\tau
({\bf P}_h \Phi^n - \Phi_h^n) \|_{L^p(\widetilde W^{-1,q})}
+\|{\bf P}_h \Phi^n - \Phi_h^n \|_{L^p(W^{1,q})}  \nonumber  \\
&\leq
C \|\Phi^n - {\bf R}_h \Phi^n \|_{L^p(W^{1,q})}  
+Ch^{-1}\|{\bf P}_h\Phi^0-\Phi^0_h\|_{L^q} \, . 
\label{lemma-5}
\end{align}
}
\end{lemma}
The estimates \eqref{lemma-3} and \eqref{lemma-4} can be found in \cite[(1.18)]{LS3} and \cite[(2.4)]{LS3}, respectively, and \eqref{lemma-5} can be proved by using \cite[(2.3)]{LS3}. 

In addition, for the elliptic boundary value problem
\begin{align}\label{Neumann-Elliptic}
\left\{\begin{aligned}
& \nabla \cdot (a \nabla u) = f + \nabla \cdot {\bf g}  &&\mbox{in}\,\,\, 
\Omega,\\
& a \nabla u \cdot {\bf n} = {\bf g} \cdot {\bf n}   
&&\mbox{on}\,\,\,\partial\Omega,
\end{aligned}\right.
\end{align}
with the constraint $\int_\Omega u\d x=0$, 
the following $W^{2,q}$ and $C^{2,\alpha}$ estimates are consequences  
of \cite[Theorem 2.4.2.7]{Grisvard} and \cite[Theorem 4.40 and  
Corollary 4.41]{Lieberman2013}.

\begin{lemma}\label{lemma3.3}
{\it
Assume that ${\bf g}=0$, $f\in L^q$ with $q \in [2, \infty)$  
and $\int_\Omega f\d x=0$, and the matrix  
$a=(a_{ij})_{d\times d}$ satisfies the ellipticity condition 
\eqref{ellipticity}. 

(1) If $a_{ij}\in W^{1,\infty}$, then \eqref{Neumann-Elliptic}  
has a unique solution $u\in W^{2,q}$ satisfying
\begin{align}\label{W2q-Elliptic}
\|u\|_{W^{2,q}}\le C_q\|f\|_{L^q},  
\end{align}
where the constant $C_q$ may depend on  
$\sum_{i,j=1}^d\|a_{ij}\|_{W^{1,\infty}}$.

(2) If $a_{ij}\in C^{1,\alpha}$, then \eqref{Neumann-Elliptic} 
has a unique solution $u\in C^{2,\alpha}$ satisfying
\begin{align}\label{C2alpha-Elliptic}
\|u\|_{C^{2,\alpha}}\le C \|f\|_{C^{\alpha}}  ,
\end{align}
where the constant $C$ may depend on $\alpha$ and 
$\sum_{i,j=1}^d\|a_{ij}\|_{C^{1,\alpha}}$.
}
\end{lemma}

Moreover, we need the following $C^{1,\alpha}$ estimate, 
which is a consequence of the steady-state case of the estimate in \cite[Theorem 4.30]{Lieberman}. 

\begin{lemma}\label{C1alpha-Elliptic}
{\it
Assume that $f \in L^\infty$, ${\bf g} \in C^{\alpha}$ for a given 
$\alpha \in (0,1)$, and $a_{ij} \in C^{\alpha}$ 
satisfies  the ellipticity condition \eqref{ellipticity}.  
Then the solution of \eqref{Neumann-Elliptic} satisfies
\begin{align}
\|u\|_{C^{1,\alpha}}\le C(\|f\|_{L^\infty}+\|{\bf g}\|_{C^{\alpha}}) , 
\end{align}
where the constant $C$ may depend on $\alpha$ and 
$\sum_{i,j=1}^d\|a_{ij}\|_{C^{\alpha}}$.
}
\end{lemma}

A $W^{1,q}$  
estimate of the corresponding finite element solution 
is given in the following lemma (a consequence of  \cite[Corollary A.6]{Gei2}).

\begin{lemma}[$W^{1,q}$ estimate of elliptic finite element  
equations]\label{Lemma:W1q-FEM}
{\it
Let $r\ge 1$, $q\in[2,\infty)$,  and ${\bf g}\in (L^q)^d$.  
If the matrix $a=(a_{ij})_{d\times d}\in W^{1,\infty}$ satisfies 
the ellipticity condition \eqref{ellipticity}, then the finite  
element system 
\begin{align}\label{FEM-Shr}
\big(a\nabla u_h ,\nabla v_h\big) =({\bf g},\nabla v_h), 
\quad\forall\, v_h\in \mathring S_h^r,
\end{align}
has a unique solution $u_h\in \mathring S_h^r$, satisfying
\begin{align}\label{W1q-FEM-est}
\|u_h\|_{W^{1,q}}\le C_q\|{\bf g}\|_{L^q} ,
\end{align}
where $C_q$ may depend on $\sum_{i,j=1}^d\|a_{ij}\|_{W^{1,\infty}}$.
}
\end{lemma}

The following discrete version of inhomogeneous Sobolev embedding (as a consequence of \cite[Proposition 1.2.10]{Lunardi95}) establishes a connection between Lemmas \ref{lemma1}-\ref{lemma2} and the $L^\infty$ boundedness of numerical solutions.

\begin{lemma}[Discrete inhomogeneous Sobolev embedding]\label{discr-embed}
{\it
Let $p,q\in(1,\infty)$ satisfy $2/p+d/q<1$, and let $\phi^n\in W^{1,q}$, $n=0,1,2,\dots$, be a sequence of functions such that $\phi^0=0$. Then for $\alpha\in(0,1-2/p-d/q)$ there holds
\begin{align}\label{ineq-discr-embed}
&\|\phi^n\|_{L^\infty(L^\infty)}+\|\phi^n\|_{L^\infty(C^{\alpha})} \le C(\|D_\tau\phi^n\|_{L^p(\widetilde W^{-1,q})} + \|\phi^n\|_{L^p(W^{1,q})}) ,\\
&\|\phi^n\|_{L^\infty(W^{1,\infty})}+\|\phi^n\|_{L^\infty(C^{1,\alpha})} \le C(\|D_\tau\phi^n\|_{L^p(L^q)} + \|\phi^n\|_{L^p(W^{2,q})}) ,
\label{ineq-discr-embed2}
\end{align}
where the constant $C$ is independent of $n\ge 1$.
}
\end{lemma}

The following lemma is an extension of the generalized  
Gr\"onwall's inequality \cite[Lemma 3.2]{LS2} to the  
time-discrete setting.
\begin{lemma}\label{Lemma:GronW}
{\it
Let $1<p<\infty$ and let $Y^n\ge 0$, $n=0,1,\dots,N$, be a  
sequence of numbers such that
\begin{align}\label{Gronw-cond}
\mbox{$\left(\tau\sum_{n=k+1}^m|Y^n|^p\right)^{\frac1p}
\leq \alpha \left(Y^{k}+\tau\sum_{n=k+1}^m Y^n \right)+\beta,$}
\qquad\forall\, 0 \le k<m \le N,
\end{align}
for some positive constants $ \alpha$ and $\beta$. Then there  
exists $\tau_p$ such that for $\tau\le \tau_p$,  
\begin{align}\label{Gronw-concl}
\mbox{$\left(\tau\sum_{n=0}^N|Y^n|^p\right)^{\frac1p}$}  
\leq C_{T,\alpha,p} (Y^0+\beta) ,
\end{align}
where the constants $\tau_p$ and $C_{T,\alpha,p}$ are independent of $\tau$, $\beta$ and the sequence $Y^n$, $n=0,1,\dots,N$.
}
\end{lemma}

Besides the lemmas above, the following interpolation inequality 
will be frequently used:
\begin{align}
\|v\|_{L^s} 
\le C_\epsilon\|v\|_{L^{s_1}} 
+\epsilon\|v\|_{L^{s_2}} , \quad\forall\, s\in(s_1,s_2) ,
\end{align}
where $\epsilon>0$ can be arbitrarily small at the expense of enlarging 
the constant $C_\epsilon$.
Since $W^{1,q}\hookrightarrow L^\infty$, it follows that 
\begin{align} \label{interp-ineq-}
\|v\|_{L^s} 
\le C\epsilon\|v\|_{L^2} 
+\epsilon\|v\|_{L^\infty}
\le C_\epsilon\|v\|_{L^2} 
+\epsilon\|v\|_{W^{1,q}} , \quad\forall\, s\in(2,\infty) .
\end{align}

\mbox{}
\vskip0.1in 
%
%
%

\section{$L^p$ estimates for a time-discrete system}
\label{sec:lptime}
\setcounter{equation}{0}
We define a time-discrete system corresponding to \eqref{e1}-\eqref{b2} by
\begin{align}
&\mbox{$
-\nabla\cdot\left({\frac{k(x)}{\mu(\mathcal{C}^{n-1})}}\nabla P^{n-1}\right)=q_{I}^{n-1}-q_{P}^{n-1},
\quad n=1,\dots,N+1,
$}
\label{tdis2} \\[10pt] 
&\,\, \gamma D_{\tau}\mathcal{C}^{n}-\nabla\cdot(D({\bf {U}}^{n-1})\nabla \mathcal{C}^{n})
+\mathcal{C}^n \nonumber\\
&
\mbox{$
=\hat{c}q_{I}^n+\bigg(1-\frac{1}{2}(q^n_I+q^n_P)\bigg)\mathcal{C}^{n}-\frac{1}{2}{\bf{U}}^{n-1}\cdot\nabla \mathcal{C}^{n}-\frac{1}{2}\nabla\cdot({\bf U}^{n-1}\mathcal{C}^n), 
\quad  n=1,\dots,N,
$}
\label{tdis1}
\end{align}
with the boundary and initial conditions
\begin{align}
&\mbox{$
D({\bf {U}}^{n-1})\nabla \mathcal{C}^{n}\cdot{\bf {n}}=0,
\quad \frac{k(x)}{\mu({\mathcal{C}^{n-1}})} \nabla P^{n-1}\cdot{\bf {n}}=0,$} && \text{for $x\in \partial\Omega$},
\label{tdis3}\\
&\mathcal{C}^0=c_0(x), && \text{for $x\in\Omega$},\label{tdis4}
\end{align}
where
\begin{align}
\mbox{$ {\bf {U}}^{n-1}=-\frac{k(x)}{\mu(\mathcal{C}^{n-1})}\nabla P^{n-1}, $}
\label{TDIS3}
\end{align}
and the condition $\int_{\Omega}P^{n-1}dx=0$ is enforced for the uniqueness of the solution of \eqref{tdis2}.

The fully discrete system \eqref{fdis2}-\eqref{TDIS4} can be viewed as the spatial discretization of \eqref{tdis2}-\eqref{TDIS3} by the FEM with P2 and P1 elements for $P^{n-1}$ and $\mathcal{C}^n$, 
respectively. The main result of this section is the following lemma on the $L^p$ and $L^\infty$ 
estimates for the time-discrete system \eqref{tdis2}-\eqref{TDIS3}. 
These estimates are needed for analyzing   
the fully discrete finite element solutions in the next section.

\begin{lemma}
\label{Lemma-TD}
{\it
Suppose that \eqref{e1}-\eqref{b2} has a unique solution satisfying  
\eqref{regu} for some $q\in(d,\infty)$, and let $p\in(2,\infty)$  
satisfy $2/p+d/q<1$. Then the time-discrete system  
\eqref{tdis2}-\eqref{TDIS3} has a unique solution   
$(P^{n},\mathcal{C}^{n})\in W^{2, q}\times W^{2, q}$,  
$n=0,1,\dots, N$, such that 
\begin{align}
&\|D_{\tau}\mathcal{C}^{N}\|_{L^p(L^q)}+\|\mathcal{C}^{N}\|_{L^p(W^{2, q})}  
\le C_{p,q} ,  \label{regu2} \\
&\|P^{n}\|_{W^{2,\infty}}+\|{\bf U}^{n}\|_{W^{1,\infty}} 
+\|D_\tau {\bf U}^{n}\|_{L^{\infty}}
+\|\mathcal{C}^{n}\|_{W^{1,\infty}} \leq C_{p,q} .
\label{regu3}
\end{align}
}
\end{lemma}

{\bf Proof}. 
For a given  
$\mathcal{C}^{n-1}\in W^{2,q}\hookrightarrow C^{1,\alpha}$, 
with $\alpha=1-d/q\in(0,1)$, we have  
$\frac{k}{\mu(\mathcal{C}^{n-1})}\in C^{1,\alpha}$. 
Then, by Lemma \ref{lemma3.3}, 
\eqref{tdis2} has  
a unique solution $P^{n-1}\in C^{2,\alpha}\hookrightarrow W^{2,\infty}$ 
such that 
\begin{align}
\|P^{n-1}\|_{C^{2,\alpha}}
\le C_{\|\mathcal{C}^{n-1}\|_{C^{1,\alpha}}} .
\end{align} 
where $C_{\|\mathcal{C}^{n-1}\|_{C^{1,\alpha}}}$  
is a constant depending on $\|\mathcal{C}^{n-1}\|_{C^{1,\alpha}}$. 
In view of \eqref{TDIS3},  ${\bf U}^{n-1}\in C^{1,\alpha} 
\hookrightarrow W^{1,\infty}$, $i.e.$, 
\begin{align}
\|{\bf U}^{n-1}\|_{C^{1,\alpha}}
\le C_{\|\mathcal{C}^{n-1}\|_{C^{1,\alpha}}} . 
\end{align} 
Thus by \cite[Theorem 2.4.2.7]{Grisvard}, the elliptic equation \eqref{tdis1} has a unique solution  
$\mathcal{C}^n\in W^{2,q}$, $i.e.$, 
\begin{align}
\|\mathcal{C}^n\|_{W^{2,q}}
\le C_{\|{\bf U}^{n-1}\|_{C^{1,\alpha}}}
\le C_{\|\mathcal{C}^{n-1}\|_{C^{1,\alpha}}}
\le C_{\|\mathcal{C}^{n-1}\|_{W^{2,q}}} . 
\end{align} 
This proves the existence and uniqueness of solutions 
$(P^{n},\mathcal{C}^{n})\in W^{2, q}\times W^{2, q}$, $n=0,1,\dots, N$. 
In particular, there exists an increasing function $\varphi:\R_+\rightarrow \R_+$ such that $\varphi(s)\ge s$ and 
\begin{align}\label{RegCn-Cn-1}
\|\mathcal{C}^n\|_{W^{2,q}} + \|P^{n}\|_{C^{2,\alpha}} +\|{\bf U}^{n}\|_{C^{1,\alpha}} \le \varphi(\|\mathcal{C}^{n-1}\|_{W^{2,q}}) . 
\end{align}

It remains to prove the quantitative regularity estimate  
\eqref{regu2}-\eqref{regu3}. To simplify the notations,  
we omit the dependence on $p$ and $q$ in the subscripts  
of the generic constant $C$.

We start with proving the following suboptimal $L^\infty$  
error estimate by mathematical induction:
\begin{align}
\|{\bf u}^{n}-{\bf U}^{n}\|_{L^\infty}+\|c^{n}-\mathcal{C}^{n}\|_{L^\infty} 
\leq \tau^{1/2} \, . 
\label{prima}
\end{align}
Since $c^0 - \mathcal{C}^0 =0$, the inequality above holds for $n=0$.  
We assume that \eqref{prima} holds for $0\le n\leq m-1$ and below, we prove 
that it also holds for $n=m$.

From \eqref{e2} and \eqref{tdis2}, we see that 
\begin{align}
\mbox{$
\nabla\cdot\left(\frac{k(x)}{\mu(c^{n-1})}\nabla(p^{n-1}-P^{n-1})\right)
=$}&\mbox{$\nabla\cdot\left(\left(\frac{k(x)}{\mu(c^{n-1})}-\frac{k(x)}{\mu(\mathcal{C}^{n-1})}\right)\nabla(p^{n-1}-P^{n-1})\right) 
$}
\nonumber \\
&+
\mbox{$\nabla\cdot\left(\left( \frac{k(x)}{\mu(\mathcal{C}^{n-1})}
-\frac{k(x)}{\mu(c^{n-1})}\right)
\nabla p^{n-1}\right)$}  .
\label{TperrE}
\end{align}
By the $W^{1,q}$ estimate of elliptic equations (see \cite[Theorem 1]{AQ}), we get 
\begin{align}
&\!\!\!\!\!\|p^{n-1}-P^{n-1}\|_{W^{1, q}} \nonumber \\
\le & 
\mbox{$
C\left\|\left(\frac{k(x)}{\mu(c^{n-1})}-\frac{k(x)}{\mu(\mathcal{C}^{n-1})}\right) 
\nabla(p^{n-1}-P^{n-1}) \right\|_{L^q} 
+C_q\left\|\left(\frac{k(x)}{\mu(\mathcal{C}^{n-1})}-\frac{k(x)}{\mu(c^{n-1})} \right) 
\nabla p^{n-1}\right\|_{L^q}
$} \nonumber \\
\le & C_{q}\|c^{n-1}-\mathcal{C}^{n-1}\|_{L^\infty}  \| \nabla(p^{n-1}-P^{n-1})\|_{L^q}
+C_{q}\|c^{n-1}-\mathcal{C}^{n-1}\|_{L^q}  \|\nabla p^{n-1}\|_{L^{\infty}} \nonumber \\
\le & C_{q}\tau^{\frac12} \|p^{n-1}-P^{n-1}\|_{W^{1, q}}
+C_{q}\|c^{n-1}-\mathcal{C}^{n-1}\|_{L^q}  ,\qquad \mbox{for}\,\,\, n=1,\dots,m, \nonumber
\end{align}
where we have used the induction assumption \eqref{prima} in the last inequality.
When $\tau \le \tau_1$ for some $\tau_1>0$, 
the last inequality further implies 
\begin{align}\label{pw1q}
\|p^{n-1}-P^{n-1}\|_{W^{1, q}} \le  C_{q}\|c^{n-1}-\mathcal{C}^{n-1}\|_{L^q}   ,\qquad \mbox{for}\,\,\, n=1,\dots,m .
\end{align}
By using \eqref{e2} and \eqref{TDIS3}, we have
\begin{align} \label{uTerr}
\|{\bf {u}}^{n}-{\bf {U}}^{n}\|_{L^s}
=&
\mbox{$
\left\|-\left(\frac{k(x)}{\mu(c^{n})}-\frac{k(x)}{\mu(\mathcal{C}^{n})}\right)
\nabla p^{n}-\frac{k(x)}{\mu(\mathcal{C}^{n})}
\nabla(p^{n}-P^{n})\right\|_{L^s}
$}
 \nonumber \\
\leq& C \|c^{n}-\mathcal{C}^{n}\|_{L^s} \|\nabla p^{n}\|_{L^\infty}
+C\|p^{n}-P^{n}\|_{W^{1,s}} \nn \\
\le & C\|c^{n}-\mathcal{C}^{n}\|_{L^s}
+C\|p^{n}-P^{n}\|_{W^{1,s}}
,\qquad \mbox{for}\,\,\, n=0,1,\dots,m ,
\end{align}
for any $s\in[1,\infty]$.

We rewrite \eqref{e1} into
\begin{align}
&\gamma D_{\tau}c^{n}-\nabla\cdot(D({\bf {u}}^{n-1})\nabla c^{n})+c^{n} \nonumber \\
&= \widehat{c}q_I^n+\left(1-\mbox{$\frac{1}{2}$}(q^n_I+q^n_P)\right)c^{n}-\mbox{$\frac{1}{2}$}{\bf u}^{n-1}\cdot\nabla c^{n}-\mbox{$\frac{1}{2}$}\nabla\cdot({\bf u}^{n-1}c^n) 
+E^{n}_{tr}, 
\label{rewu} 
\end{align}
where
\begin{align}
E^{n}_{tr}=
&\gamma D_{\tau}c^n-\gamma c^n_t+\nabla\cdot((D({\bf u}^{n})-D({\bf u}^{n-1}))\nabla c^n)  +({\bf u}^{n-1}-{\bf u}^n)\cdot\nabla c^n\nonumber\\
& \mbox{$ -\frac{1}{2}((q^n_I-q^n_P)-(q^{n-1}_I-q^{n-1}_P))c^n $}
\nn 
\end{align}
denotes the truncation error, satisfying the following estimate under the regularity assumption \eqref{regu}:
\begin{align} 
\|E^{n}_{tr}\|_{L^p(\widetilde W^{-1,q})}\leq C\tau.
\nn 
\end{align}
Subtracting \eqref{tdis1} from \eqref{rewu} gives 
\begin{align}
& \gamma D_{\tau}(c^{n}-\mathcal{C}^{n})-\nabla\cdot(D({\bf {u}}^{n-1})
\nabla(c^{n}-\mathcal{C}^{n}))+c^{n}-\mathcal{C}^{n} \label{eeqrt}
\\
=&\left(1-\mbox{$\frac{1}{2}$}(q_I^n+q^n_P)\right)(c^{n}-\mathcal{C}^{n})
+\mbox{$\frac{1}{2}$}(q^{n-1}_I-q^{n-1}_P)(c^{n}-\mathcal{C}^{n})-\mbox{$\frac{1}{2}$}({\bf {u}}^{n-1}
-{\bf {U}}^{n-1})\cdot\nabla c^{n}
\nn\\
& -\mbox{$\frac{1}{2}$}\nabla\cdot({\bf {U}}^{n-1}(c^{n} 
-\mathcal{C}^{n}))-\mbox{$\frac{1}{2}$}\nabla\cdot(({\bf u}^{n-1}-{\bf U}^{n-1})c^n+{\bf U}^{n-1}(c^n 
-\mathcal{C}^n))\nonumber\\
&+\nabla\cdot((D({\bf{u}}^{n-1})-D({\bf {U}}^{n-1}))\nabla (\mathcal{C}^{n}-c^n))
 +\nabla\cdot((D({\bf{u}}^{n-1})-D({\bf {U}}^{n-1}))\nabla c^n)+E^{n}_{tr}.  \nonumber
\end{align}
Applying Lemma \ref{lemma1} to the last equation yields, for $p \in (2, \infty)$ and $n=1,\dots,m$,
\begin{align}
&\!\!\!\!\!\|D_{\tau}(c^{n}-\mathcal{C}^{n})\|_{L^p(\widetilde{W}^{-1, q})}
+\|c^{n}-\mathcal{C}^{n}\|_{L^p(W^{1, q})} \nonumber \\
\le & 
\mbox{$
C \left\|\left(1-\frac{1}{2}(q_I^n+q^n_P)\right)(c^{n}-\mathcal{C}^{n})\right\|_{L^p(L^q)}
      +C \|(q^{n-1}_I-q^{n-1}_P)(c^{n}-\mathcal{C}^{n})\|_{L^p(L^q)}
      $}
 \nonumber \\
&
+ C \|{\bf {u}}^{n-1}-{\bf {U}}^{n-1}\|_{L^p(L^q)} \|\nabla c^{n}\|_{L^{\infty}(L^{\infty})}
+ C \|{\bf {U}}^{n-1}\|_{L^{\infty}(L^{\infty})} \| c^{n}-\mathcal{C}^{n} \|_{L^p{(L^q)}}  \nonumber \\
& +C\|{\bf u}^{n-1}-{\bf U}^{n-1}\|_{L^p(L^q)}\|c^n\|_{L^{\infty}(L^{\infty})}
+C\|{\bf U}^{n-1}\|_{L^{\infty}(L^{\infty})}\|c^n-\mathcal{C}^n\|_{L^p(L^q)}\nonumber\\
&+C \|D({\bf{u}}^{n-1})-D({\bf {U}}^{n-1})\|_{L^{\infty}(L^{\infty})} \|\nabla (\mathcal{C}^{n}-c^{n})\|_{L^p{(L^q)}}  \nonumber \\
&+C \|D({\bf{u}}^{n-1})-D({\bf {U}}^{n-1})\|_{L^p(L^q)} \|\nabla c^{n}\|_{L^{\infty}(L^{\infty})}
+ C\|E^{n}_{tr}\|_{L^p{(\widetilde W^{-1,q})}}  \nonumber \\
\leq& C (\|c^{n-1}-\mathcal{C}^{n-1}\|_{L^p(L^q)}+\|c^{n}-\mathcal{C}^{n}\|_{L^p(L^q)})+C\tau^{1/2} \|c^{n}-\mathcal{C}^{n}\|_{L^p(W^{1,q})}+C\tau \nonumber \\
\leq& C \|c^{n}-\mathcal{C}^{n}\|_{L^p(L^q)} +C\tau^{1/2} \|c^{n}-\mathcal{C}^{n}\|_{L^p(W^{1,q})}+C\tau,
\label{wpesti}
\end{align}
where we have used induction assumption \eqref{prima}  to estimate $\|D({\bf{u}}^{n-1})-D({\bf {U}}^{n-1})\|_{L^{\infty}(L^{\infty})}$ and $\|{\bf {U}}^{n-1}\|_{L^{\infty}(L^{\infty})}$, and used \eqref{pw1q}-\eqref{uTerr} to estimate $\|{\bf {u}}^{n-1}-{\bf {U}}^{n-1}\|_{L^p(L^q)} $.
When $\tau \le \tau_2$ for some $\tau_2>0$, the last inequality reduces to 
\begin{align}
&\|D_{\tau}(c^{n}-\mathcal{C}^{n})\|_{L^p(\widetilde{W}^{-1, q})}
+\|c^{n}-\mathcal{C}^{n}\|_{L^p(W^{1, q})} \le C \|c^{n}-\mathcal{C}^{n}\|_{L^p(L^q)} +C\tau ,
\quad
n=1,\dots,m .
\label{wpesti22}
\end{align}
By Lemma \ref{discr-embed},  
\begin{align}
\begin{aligned}
\|c^{n}-\mathcal{C}^{n}\|_{L^{\infty}({L}^\infty)}
&\le C(\|D_{\tau}(c^{n}-\mathcal{C}^{n})\|_{L^p(\widetilde{W}^{-1, q})}
+\|c^{n}-\mathcal{C}^{n}\|_{L^p(W^{1, q})} )  \\
&\le C  \|c^{n}-\mathcal{C}^{n}\|_{L^p(L^q)} +C\tau \qquad\mbox{[\eqref{wpesti22} is used here]} \\
&\le C  \|c^{n}-\mathcal{C}^{n}\|_{L^p(L^\infty)} +C\tau   \\
&\le \mbox{$\frac{1}{2}$} \|c^{n}-\mathcal{C}^{n}\|_{L^\infty(L^\infty)} 
+ C_{p,q} \|c^{n}-\mathcal{C}^{n}\|_{L^1(L^\infty)} +C\tau ,
\quad n=1,\dots,m,
\end{aligned}
\nn 
\end{align}
which further implies (through applying Gronwall's inequality, $i.e.$, Lemma \ref{Lemma:GronW})
\begin{align}\label{cm-Cm-Linfty} 
\|c^{m}-\mathcal{C}^{m}\|_{L^{\infty}({L}^\infty)} \le C\tau .
\end{align}
Substituting the inequality above into \eqref{wpesti22}, we have
\begin{align}\label{cm-Cm-LpW1q1}
&\|D_{\tau}(c^{m}-\mathcal{C}^{m})\|_{L^p(\widetilde{W}^{-1, q})}
+\|c^{m}-\mathcal{C}^{m}\|_{L^p(W^{1, q})}
\le C \tau  
\end{align}
which with \eqref{ineq-discr-embed} shows 
\begin{align}\label{cm-Cm-LinfCalpha}
&\|c^{m}-\mathcal{C}^{m}\|_{L^\infty(C^{\alpha})}
\le C(\|D_{\tau}(c^{m}-\mathcal{C}^{m})\|_{L^p(\widetilde{W}^{-1, q})}
+\|c^{m}-\mathcal{C}^{m}\|_{L^p(W^{1, q})} )
\le C \tau .
\end{align}
By using an inverse inequality in time, \eqref{cm-Cm-LpW1q1} implies
\begin{align}\label{cm-Cm-W1q1}
\|c^{m}-\mathcal{C}^{m}\|_{L^\infty(W^{1, q})} \le C \tau^{1-1/p} .
\end{align}


Moreover, applying \eqref{W2q-Elliptic} to \eqref{TperrE} leads to, for $n=0,1,\dots,m$,  
\begin{align} \label{eptau-W2q}
&\|p^{n}-P^{n}\|_{W^{2,q}}  \\
 &\le 
\mbox{$C \left\|\nabla\cdot\left(\left(\frac{k(x)}{\mu(c^{n})}-\frac{k(x)}{\mu(\mathcal{C}^{n})}\right)
\nabla(p^{n}-P^{n})\right)\right\|_{L^q} +
C\left\|\nabla\cdot\left(\left( \frac{k(x)}{\mu(\mathcal{C}^{n})}
-\frac{k(x)}{\mu(c^{n})}\right)
\nabla p^{n}\right) \right\|_{L^q} 
$}
\nonumber \\
&\le 
C 
(\|c^{n}-\mathcal{C}^{n}\|_{L^\infty}\|p^{n}-P^{n}\|_{W^{2,q}}
+\|c^{n}-\mathcal{C}^{n}\|_{W^{1,q}}\|p^{n}-P^{n}\|_{W^{1,\infty}}) \nonumber \\
&\quad +C 
(\|c^{n}-\mathcal{C}^{n}\|_{L^\infty}\|p^{n}\|_{W^{2,q}}
+\|c^{n}-\mathcal{C}^{n}\|_{W^{1,q}}\|p^{n}\|_{W^{1,\infty}})  \nonumber \\
&\le 
C (\tau \|p^{n}-P^{n}\|_{W^{2,q}} + \tau^{1-1/p} \|p^{n}-P^{n}\|_{W^{2,q}}) +C(\tau +\|c^{n}-\mathcal{C}^{n}\|_{W^{1,q}}) ,
\nn 
\end{align}
where we used \eqref{cm-Cm-Linfty}-\eqref{cm-Cm-W1q1} in deriving the last inequality. 
When $\tau \le \tau_3$ for some $\tau_3>0$, we see that 
\begin{align}
\|p^{n}-P^{n}\|_{W^{2,q}}  \le &
C(\tau +\|c^{n}-\mathcal{C}^{n}\|_{W^{1,q}}) ,\quad n=0,1,\dots,m.
\nn 
\end{align}
By noting  
\eqref{cm-Cm-LpW1q1} 
and the Sobolev embedding $W^{2,q}\hookrightarrow W^{1,\infty}$ for $q>d$, we obtain
\begin{align}
\|p^{m}-P^{m}\|_{L^p(W^{1,\infty})} \le C \|p^{m}-P^{m}\|_{L^p(W^{2,q})} \le &
C ( \tau +\|c^{m}-\mathcal{C}^{m}\|_{L^p(W^{1,q})} )\le C \tau 
\end{align}
which, together with an  inverse inequality in time, leads to 
\begin{align}\label{pm-Pm-LinftyW2q}
\|p^{m}-P^{m}\|_{L^\infty(W^{1,\infty})}
+\|p^{m}-P^{m}\|_{L^\infty(W^{2,q})}  \le C\tau^{1-1/p} .
\end{align}
By taking $s=\infty$ in \eqref{uTerr} and using \eqref{cm-Cm-Linfty}, we get
\begin{align}\label{um-Um-Linfty}
\|{\bf {u}}^{m}-{\bf {U}}^{m}\|_{L^\infty(L^\infty)}  
\le C\tau^{1-1/p}  . 
\end{align}
Since $\frac{2}{p}+\frac{d}{q}<1$ implies $p>2$ and therefore $C\tau^{1-1/p}\le \tau^{1/2}$ for sufficiently small stepsize $\tau$, by combining above result and \eqref{cm-Cm-Linfty}, the mathematical induction on \eqref{prima} is closed as $\tau\le\tau_4$ for some  $\tau_4>0$. 
Consequently, the estimates \eqref{cm-Cm-Linfty}, \eqref{cm-Cm-W1q1},  
\eqref{pm-Pm-LinftyW2q} 
and \eqref{um-Um-Linfty} hold for $m=N$. 
When $\displaystyle \tau \le \min_{1 \le j \le 4} \tau_j$, 
we have the following estimates:
\begin{align}\label{cm-Cm-LpW1q1-2}
\|\mathcal{C}^{n}\|_{L^\infty}+\|\mathcal{C}^{n}\|_{W^{1, q}}
+ \|P^{n}\|_{W^{1,\infty}} +\|P^{n}\|_{W^{2,q}}
+\|{\bf {U}}^{n}\|_{L^\infty} +\|D_{\tau} C^n \|_{L^\infty}  \le C .
\end{align}
From \eqref{TDIS3} we further see that
\begin{align}\label{um-Um-Linfty-2}
\|{\bf {U}}^{n}\|_{W^{1,q}}
\le C(\|P^{n}\|_{W^{2,q}}+\|\mathcal{C}^{n}\|_{W^{1,q}}  \|P^{n} \|_{W^{1,\infty}})
\le C, \quad n=0,1,...,N \, .
\end{align}

Now we are ready to prove \eqref{regu2}-\eqref{regu3}. 
To prove \eqref{regu2}, we rewrite \eqref{tdis1} into
\begin{align}
&\gamma D_{\tau}\mathcal{C}^{n}-\nabla\cdot(D({\bf {u}}^{n-1})\nabla \mathcal{C}^{n})
+\mathcal{C}^n \nn \\
&=\hat{c}q_{I}^n+\left(1-\mbox{$\frac{1}{2}$}(q^n_I+q^n_P)\right)\mathcal{C}^{n}-\mbox{$\frac{1}{2}$}{\bf{U}}^{n-1}\cdot\nabla \mathcal{C}^{n}-\mbox{$\frac{1}{2}$}\nabla\cdot({\bf U}^{n-1}\mathcal{C}^n) \nn \\
&\quad +\nabla\cdot((D({\bf U^{n-1}})-D({\bf u}^{n-1}))\nabla\mathcal{C}^n)    \nn
\end{align}
and by Lemma \ref{lemma1}, we obtain
{\small
\begin{align}\label{Pre-LpqCN}
&\!\!\!\!\! \|D_{\tau}\mathcal{C}^{N}\|_{L^p(L^q)}+\|\mathcal{C}^{N}\|_{L^p(W^{2, q})} \\
\le & C\|\hat{c}q_I^N\|_{L^p(L^q)}+C\left\|(1-\mbox{$\frac{1}{2}$}(q^N_I+q^N_P))\mathcal{C}^{N}\right\|_{L^p(L^q)}
+C\|{\bf {U}}^{N-1}\|_{L^{\infty}(L^{\infty})}\|\nabla \mathcal{C}^{N}\|_{L^p(L^q)}
\nn\\
&+C\|\nabla\cdot({\bf U}^{N-1}\mathcal{C}^N)\|_{L^p(L^q)}    +C\|\nabla\cdot((D({\bf {U}}^{N-1})-D({\bf {u}}^{N-1}))
\nabla \mathcal{C}^{N})\|_{L^p(L^q)}
\nonumber\\
\le &C+
C\|\nabla\cdot((D({\bf {U}}^{N-1})-D({\bf {u}}^{N-1}))
\nabla \mathcal{C}^{N})\|_{L^p(L^q)})  
&&\hskip-4cm \mbox{(use \eqref{cm-Cm-LpW1q1-2}-\eqref{um-Um-Linfty-2})}  \nonumber \\
\le &C+ C (\|\nabla {\bf U}^{N-1}\|_{L^\infty(L^q)}+\|\nabla {\bf u}^{N-1}\|_{L^\infty(L^q)}) \|\nabla \mathcal{C}^{N}\|_{L^p(L^{\infty})} \nonumber \\
&+C\|D({\bf {U}}^{N-1})-D({\bf {u}}^{N-1})\|_{L^{\infty}(L^{\infty})} \|\mathcal{C}^{N}\|_{L^p(W^{2,q})} 
\nn\\
\le &C+C\|\mathcal{C}^{N}\|_{L^p(W^{1,\infty})}+C\tau^{\frac12}\|\mathcal{C}^{N}\|_{L^{p}(W^{2, q})} . && \hskip-4cm
\mbox{(use \eqref{prima} and \eqref{um-Um-Linfty-2})}  \nonumber
\end{align}
}
By noting $\|\mathcal{C}^{n}\|_{W^{1,\infty}} 
\le \frac{1}{2} \|\mathcal{C}^{n}\|_{W^{2,q}}  
+C\|\mathcal{C}^{n}\|_{W^{1,q}}$ and \eqref{cm-Cm-LpW1q1-2},  
when $\tau \le \tau_5$ for some $\tau_5>0$, 
\eqref{Pre-LpqCN} reduces to
\begin{align}\label{LpW2qCn-1}
&\!\!\!\!\! \|D_{\tau}\mathcal{C}^{N}\|_{L^p(L^q)} 
+\|\mathcal{C}^{N}\|_{L^p(W^{2, q})}
\le
C+C\|\mathcal{C}^{N}\|_{L^p(W^{1,q})}
\le C \, . 
\end{align}
\eqref{regu2} is obtained.

To prove \eqref{regu3}, we use \eqref{LpW2qCn-1} and Lemma \ref{discr-embed}, which imply
\begin{align} \label{CnC1alpha}
\|\mathcal{C}^{N}\|_{L^\infty(C^{1,\alpha})}
\le C(\|D_{\tau}\mathcal{C}^{N}\|_{L^p(L^q)}+\|\mathcal{C}^{N}\|_{L^p(W^{2, q})} )
\le C . 
\end{align}
With the regularity estimate above, applying   [Lemma \ref{lemma3.3}, \eqref{C2alpha-Elliptic}]  to \eqref{tdis2} yields
\begin{align}\label{PnC2alpha}
\|P^{n}\|_{C^{2,\alpha}}
\le C\|q_I^n-q_P^n\|_{C^\alpha} \le C,\quad n=0,1,\dots,N ,
\end{align}
and substituting \eqref{CnC1alpha}-\eqref{PnC2alpha} into \eqref{TDIS3} gives 
\begin{align}\label{UnC1alpha}
\|{\bf U}^{n}\|_{C^{1,\alpha}}  \le C,\quad n=0,1,\dots,N.
\end{align}
Again, applying the backward difference operator $D_\tau$ to \eqref{tdis2} yields
\begin{align}
&
\mbox{$
-\nabla\cdot\left({\frac{k(x)}{\mu(\mathcal{C}^{n})}}\nabla D_\tau P^{n}\right)
-\nabla\cdot\left(D_{\tau}\left({\frac{k(x)}{\mu(\mathcal{C}^{n})}}\right)\nabla P^{n-1}\right)
=D_\tau q_{I}^{n}-D_\tau q_{P}^{n} .
$}
\end{align}
By Lemma \ref{C1alpha-Elliptic}, 
\begin{align}\label{DtauPnLinfty}
\|D_\tau P^{n}\|_{C^{1,\alpha}}
&\le 
\mbox{$ 
\frac{C}{\tau}\left\|\left({\frac{k(x)}{\mu(\mathcal{C}^{n})}}-{\frac{k(x)}{\mu(\mathcal{C}^{n-1})}}\right)\nabla P^{n-1}\right\|_{C^{\alpha}}
+C\|D_\tau q_{I}^{n}-D_\tau q_{P}^{n}\|_{L^\infty}
$}
 \nonumber \\
&\le \mbox{$ 
\frac{C}{\tau}$} \|\mathcal{C}^{n}-\mathcal{C}^{n-1}\|_{C^{\alpha}}\|\nabla P^{n-1}\|_{C^{\alpha}}
+C(\|\partial_t q_{I}\|_{L^\infty(0,T;L^\infty)}+\|\partial_t q_{P}\|_{L^\infty(0,T;L^\infty)}) \nonumber \\
&\le \mbox{$ 
\frac{C}{\tau} \|c^{n}-c^{n-1}\|_{C^{\alpha}}+\frac{C}{\tau}\|\mathcal{C}^{n}-c^{n}\|_{C^{\alpha}} +\frac{C}{\tau}\|\mathcal{C}^{n-1}-c^{n-1}\|_{C^{\alpha}}+C
$} \nonumber \\
&\le C, \qquad n=1,\dots,N,
\end{align}
where we have used \eqref{cm-Cm-LinfCalpha} in the last inequality. Finally,  
from \eqref{TDIS3} we see that 
\begin{align}
\mbox{$
\| D_\tau{\bf {U}}^{n} \|_{L^\infty} \le 
C\left( \| \nabla D_\tau P^{n}\|_{L^\infty} + 
\left\| D_{\tau}\left({\frac{k(x)}{\mu(\mathcal{C}^{n})}}\right) \right \|_{L^\infty}\right) \le C ,
\quad n=1,\dots,N,
$}
\end{align}
and \eqref{regu3} follows immediately.
This proves Lemma \ref{Lemma-TD} in the case  
$\displaystyle \tau \le \tau_{p,q}^*:=\min_{1 \le j \le 4} \tau_j$.

If $\displaystyle \tau \ge \tau_{p,q}^*$, 
$N=T/\tau\le T/ \tau_{p,q}^*\le C$, and therefore, \eqref{RegCn-Cn-1} implies 
\begin{align}
\|\mathcal{C}^n\|_{W^{2,q}} +\|P^{n}\|_{C^{2,\alpha}}   
+\|{\bf U}^{n}\|_{C^{1,\alpha}} 
\le \varphi^{(n)}(\|\mathcal{C}^{0}\|_{W^{2,q}})
\le  \varphi^{(T/ \tau_{p,q}^*)}(\|\mathcal{C}^{0}\|_{W^{2,q}}) \le C ,
\end{align} 
where $ \varphi^{(n)}:= \varphi^{(n-1)}\circ \varphi$. 
This proves Lemma \ref{Lemma-TD} in the case  
$\displaystyle \tau \ge \tau_{p,q}^*$. 
\endproof

\section{The proof of Theorem \ref{MainTHM}}
\label{sec:lpfully}
\setcounter{equation}{0}

Before proving Theorem \ref{MainTHM}, 
we show the boundedness of the numerical solutions based on the uniform 
regularity estimates given in Lemma \ref{Lemma-TD}  
for the time-discrete system \eqref{tdis2}-\eqref{TDIS3}.

\subsection{Boundedness of the numerical solutions}$\,$

\begin{lemma}\label{Lemma:BD}
Under the assumption of Theorem \ref{MainTHM}, there exist positive constants  
$\tau_q$ and $h_q$ such that for $\tau\le \tau_q$ and $h\leq h_q$ the finite   
element system \eqref{fdis2}-\eqref{TDIS4} has a unique solution
$(P^{n}_h, \mathcal{C}^{n}_h)$, $n=0,1,..., N$, satisfying the following  
estimates:
\begin{align}
\|\mathcal{C}^n_h\|_{W^{1,\infty}} +\|{\bf U}^n_h \|_{L^\infty} \le C .
\label{cCerr}
\end{align}
\end{lemma}

{\bf Proof.} Since both coefficient matrices of the linear systems 
\eqref{fdis2} and \eqref{fdis1} are positive definite (possibly non-symmetric), 
it follows that the linear system \eqref{fdis2}-\eqref{fdis1} has a  
unique solution. 

Next, we prove a primary estimate 
\begin{align} 
\label{priest}
\|{\bf{P}}_h\mathcal{C}^n-\mathcal{C}^n_h\|_{L^\infty}\leq h^{\frac12}  , 
\quad n=0,\dots,m-1,
\end{align}
by mathematical induction. 
For the given $q>d$, we choose a fixed $p\in(2,\infty)$ satisfying  
$2/p+d/q<1$, and omit the dependence on $p$ and $q$ in the subscripts  
of generic constants below.

Since
$\|{\bf{P}}_h\mathcal{C}^0-\mathcal{C}^0_h\|_{L^{\infty}}
= \|{\bf{P}}_h c_0 - \Pi_h c_0\|_{L^{\infty}}
\le C h\|c_0\|_{W^{1,{\infty}}}$,  
(\ref{priest}) holds for $m=1$ when $h \le h_1$ for some $h_1>0$. 
Therefore, we can assume that it holds for some positive integer $m$.

From \eqref{tdis2}, we see that 
\begin{align}
&\mbox{$ 
 \left({\frac{k(x)}{\mu(\mathcal{C}^{n-1})}}\nabla P^{n-1}, \nabla v_h\right)
=(q_I^{n-1}-q_P^{n-1}, v_h),
\qquad\qquad\qquad \forall\, v_h\in  \mathring S^{2}_h .
$}
\end{align}
and therefore, subtracting the equation above from \eqref{fdis2} yields 
\begin{align}
\mbox{$ \nabla\cdot\left(\frac{k(x)}{\mu(\mathcal{C}^{n-1})} 
\nabla(\overline{\bf P}_h P^{n-1}-P_h^{n-1})\right)
$}
=&
\mbox{$ 
\nabla\cdot\left(\left(\frac{k(x)}{\mu(\mathcal{C}^{n-1})} 
-\frac{k(x)}{\mu(\mathcal{C}_h^{n-1})}\right)\nabla(P^{n-1}-P_h^{n-1})\right)
$}  \nonumber \\
&\mbox{$ 
+
\nabla\cdot\left(\left( \frac{k(x)}{\mu(\mathcal{C}_h^{n-1})}
-\frac{k(x)}{\mu(\mathcal{C}^{n-1})}\right)
\nabla P^{n-1}\right)
$}  \nonumber \\
&\mbox{$ 
+ \nabla\cdot \left(\frac{k(x)}{\mu(\mathcal{C}^{n-1})} 
\nabla(\overline{\bf P}_h P^{n-1}-P^{n-1})\right)  
$} .
\label{TperrE-FE}
\end{align}
Since $\|\frac{k(x)}{\mu(\mathcal{C}^{n-1})}\|_{W^{1,\infty}}\le C$ (as a consequence of [Lemma \ref{Lemma-TD}, \eqref{regu3}]), by the $W^{1,s}$ estimate of elliptic finite element system (Lemma \ref{Lemma:W1q-FEM}), we have 
\begin{align} \label{pnw1q}
&\|P^{n-1}-P^{n-1}_h\|_{W^{1, s}} \\
&\le 
\mbox{$
C\left\|\left(\frac{k(x)}{\mu(\mathcal{C}^{n-1})}-\frac{k(x)}{\mu(\mathcal{C}_h^{n-1})}\right)\nabla(P^{n-1}-P_h^{n-1})\right\|_{L^s}
     +C\left\|\left( \frac{k(x)}{\mu(\mathcal{C}_h^{n-1})} -\frac{k(x)}{\mu(\mathcal{C}^{n-1})}\right)  \nabla P^{n-1} \right\|_{L^s}
     $}
\nonumber\\
&\quad +C\|\overline{\bf P}_h P^{n-1}-P^{n-1}\|_{W^{1, s}} \nonumber\\
&\le C\|\mathcal{C}^{n-1}-\mathcal{C}_h^{n-1}\|_{L^\infty} \|P^{n-1}-P^{n-1}_h\|_{W^{1, s}}
     +C\|\mathcal{C}^{n-1}-\mathcal{C}_h^{n-1}\|_{L^s} \|P^{n-1}\|_{W^{1, \infty}} +Ch\|P^{n-1}\|_{W^{2, s}}
\nonumber\\
&\le   Ch^{\frac12} \|P^{n-1}-P^{n-1}_h\|_{W^{1, s}}
     +C\|\mathcal{C}^{n-1}-\mathcal{C}_h^{n-1}\|_{L^s} +Ch ,
\qquad n=1,\dots,m ,
\quad\forall\, s\in(1,\infty),
\nn
\end{align}
where we have used the induction assumption \eqref{priest} to estimate  
$\|\mathcal{C}^{n-1}-\mathcal{C}_h^{n-1}\|_{L^\infty}$, and  
Lemma \ref{Lemma-TD} to estimate $\|P^{n-1}\|_{W^{1, \infty}}$  
and $\|P^{n-1}\|_{W^{2, s}}$.
Choosing $s=4d$ in the last equation, we can see that 
when $h \le h_2$ for some $h_2>0$, 
\begin{align}\label{pnw1q22}
&\|P^{n-1}-P^{n-1}_h\|_{W^{1,4d}} \le C\|\mathcal{C}^{n-1} 
-\mathcal{C}_h^{n-1}\|_{L^{4d}} +Ch ,
\qquad n=1,\dots,m . 
\end{align}
By an inverse inequality, 
\begin{align}\label{P-Ph-W1infty}
\|P^{n-1}-P^{n-1}_h\|_{W^{1, \infty}}
&\le  \|P^{n-1}-\overline{\bf P}_h P^{n-1}\|_{W^{1, \infty}} + \|\overline{\bf P}_h P^{n-1}-P^{n-1}_h\|_{W^{1, \infty}} \nonumber \\
&\le  Ch\|P^{n-1}\|_{W^{2, \infty}} + Ch^{-\frac14}\|\overline{\bf P}_h P^{n-1}-P^{n-1}_h\|_{W^{1, {4d}}}
\nonumber \\
&\le  Ch\|P^{n-1}\|_{W^{2, \infty}} + Ch^{-\frac14}(\|\mathcal{C}^{n-1}-\mathcal{C}_h^{n-1}\|_{L^{4d}} +h)
&&\mbox{(use \eqref{pnw1q22} here)} 
\nonumber \\
&\le  Ch + Ch^{-\frac14}(h^{\frac12} +h) 
&&\mbox{(use \eqref{priest} here)} \nonumber \\
&\le  Ch^{\frac14} ,
\qquad n=1,\dots,m,
\end{align}
where we have used Lemma \ref{Lemma-TD}  and the induction assumption \eqref{priest}. 
Moreover, subtracting \eqref{TDIS3} from \eqref{TDIS4}  and using Lemma 3.1 and Lemma 4.1,  
we derive
\begin{align} \label{uerr}
&\|{\bf{U}}^{n-1}-{\bf{U}}^{n-1}_h\|_{L^s} \nonumber \\
&\le 
\mbox{$
\left\|\frac{k(x)}{\mu(\mathcal{C}^{n-1}_h)}\nabla(P^{n-1}_h-P^{n-1})
+\left(\frac{k(x)}{\mu({\mathcal{C}^{n-1}_h})}-\frac{k(x)}{\mu({\mathcal{C}^{n-1}})}\right)
\nabla  P^{n-1}\right\|_{L^s} 
$}
\nn\\
&\le
C\|P^{n-1}-P^{n-1}_h\|_{W^{1,s}} + C\|\mathcal{C}^{n-1}-\mathcal{C}^{n-1}_h\|_{L^{s}}
\|P^{n-1}\|_{W^{1, \infty}}
\nn\\
&\le
C\|P^{n-1}-P^{n-1}_h\|_{W^{1,s}}+C\|\mathcal{C}^{n-1}-\mathcal{C}^{n-1}_h\|_{L^{s}} ,
\qquad n=1,\dots,m,
\quad\forall\, s\in[1,\infty] .
\end{align}
Setting $s=\infty$ in the inequality above and 
using \eqref{P-Ph-W1infty} and the induction assumption \eqref{priest},  
we obtain 
\begin{align} \label{U-Uh_LinftyLinf}
&\|{\bf{U}}^{n-1}-{\bf{U}}^{n-1}_h\|_{L^\infty} \nonumber \\
&\le
C\|P^{n-1}-P^{n-1}_h\|_{W^{1,\infty}}+C\|\mathcal{C}^{n-1}-\mathcal{C}^{n-1}_h\|_{L^\infty}
\le  Ch^{\frac14},
\quad n=1,\dots,m \, . 
\end{align}
Similarly, choosing $s=q$ in \eqref{pnw1q} and \eqref{uerr}, we have 
\begin{align}\label{U-Uh_LpLq}
&\|{\bf{U}}^{n-1}-{\bf{U}}^{n-1}_h\|_{L^q} + \|P^{n-1}-P^{n-1}_h\|_{W^{1,q}} \nonumber \\
&\le
C\|\mathcal{C}^{n-1}-\mathcal{C}^{n-1}_h\|_{L^q} + Ch ,
\quad n=1,\dots,m.
\end{align}


To estimate $\|\mathcal{C}^{n-1}-\mathcal{C}^{n-1}_h\|_{L^q}$, we rewrite 
the finite element system \eqref{fdis1} as
\begin{align} \label{ReWrFD-C}
& (\gamma D_{\tau}\mathcal{C}_h^{n}, w_h)
+(D({\bf {U}}^{n-1})\nabla \mathcal{C}_h^{n}, \nabla w_h) +(\mathcal{C}^n_h, w_h) \nonumber\\
&=\Big(\hat{c}q_I^n + \big(1-\mbox{$\frac{1}{2}$}(q_I^n+q_P^n)\big)\mathcal{C}^n, w_h\Big)
- \mbox{$\frac12$} ({\bf{U}}^{n-1}\cdot\nabla \mathcal{C}^n , w_h) + \mbox{$\frac{1}{2}$}({\bf{U}}^{n-1} \cdot \nabla w_h, \mathcal{C}^n )  \nonumber \\
&\quad+\big((D({\bf {U}}^{n-1})-D({\bf {U}}_h^{n-1}))\nabla \mathcal{C}_h^{n}, \nabla w_h\big)
+ \Big(\big(1-\mbox{$\frac{1}{2}$}(q_I^n+q_P^n)\big)(\mathcal{C}^n_h-\mathcal{C}^n), w_h\Big)  \nonumber \\
&\quad - \mbox{$\frac12$} (({\bf{U}}^{n-1}_h-{\bf{U}}^{n-1})\cdot\nabla \mathcal{C}^n_h, w_h)
            + \mbox{$\frac12$} ((q^{n-1}_I-q^{n-1}_P) (\mathcal{C}^n_h-\mathcal{C}^n), w_h)\nonumber \\
&\quad + \mbox{$\frac12$} (({\bf{U}}^{n-1}_h-{\bf{U}}^{n-1}) \cdot \nabla w_h, \mathcal{C}^n_h )
            + 
            ({\bf{U}}^{n-1}\cdot \nabla w_h, \mathcal{C}^n_h-\mathcal{C}^n) \, ,
\qquad\forall\, w_h\in S^1_h.
\end{align}
In view of the difference between the right-hand sides of  \eqref{tdis1} and \eqref{ReWrFD-C}, and in order to invoke Lemma 3.2, 
we define $\theta^{n}$ to be the solution of the following auxiliary time-discrete equation
 \begin{align} \label{auxeq}
&\gamma D_{\tau}\theta^{n}-\nabla\cdot(D({\bf{U}}^{n-1})\nabla\theta^{n}) +\theta^{n} \nonumber \\
&=-\nabla\cdot \big(D({\bf {U}}^{n-1})-D({\bf {U}}_h^{n-1}))\nabla \mathcal{C}_h^{n}\big)
     +\big(1-\mbox{$\frac{1}{2}$}(q_I^n+q_P^n)\big)(\mathcal{C}^n_h-\mathcal{C}^n)  \nonumber \\
&\quad\,\, - \mbox{$\frac12$} ({\bf{U}}^{n-1}_h-{\bf{U}}^{n-1})\cdot\nabla \mathcal{C}^n_h
            + \mbox{$\frac12$} (q^{n-1}_I-q^{n-1}_P) (\mathcal{C}^n_h-\mathcal{C}^n) \nonumber \\
&\quad\,\, - \mbox{$\frac12$} \nabla\cdot\big(({\bf{U}}^{n-1}_h-{\bf{U}}^{n-1}) \mathcal{C}^n_h \big)
            - \nabla\cdot\big({\bf{U}}^{n-1}(\mathcal{C}^n_h-\mathcal{C}^n)\big),
\end{align}
with the boundary and initial conditions
\begin{align*}
&-D({\bf{U}}^{n-1})\nabla \theta^{n}\cdot{\bf {n}} =-(D({\bf {U}}^{n-1})-D({\bf {U}}_h^{n-1}))\nabla \mathcal{C}^{n}_h \cdot{\bf {n}} -\mbox{$\frac{1}{2}$}({\bf{U}}_h^{n-1}-{\bf{U}}^{n-1}) \mathcal{C}^n_h \cdot{\bf n}  \nonumber \\
&\qquad\qquad\qquad\qquad\qquad\!\!\! - {\bf{U}}^{n-1}(\mathcal{C}^n_h-\mathcal{C}^n) \cdot{\bf n}
&& \mbox{on}\,\,\,\partial\Omega , \\
&\theta^0=0  &&\mbox{in}\,\,\,\Omega ,
\end{align*}
and define $\theta_h^n\in S_h^1$ to be the solution of the corresponding  
fully-discrete finite element system:
\begin{align} \label{fauxeq}
&(\gamma D_{\tau}\theta_h^{n}, w_h)
+((D({\bf{U}}^{n-1})\nabla\theta^{n}_h, \nabla w_h)
+(\theta^{n}_h, w_h)  \nonumber \\
&=\big((D({\bf {U}}^{n-1})-D({\bf {U}}_h^{n-1}))\nabla \mathcal{C}_h^{n}, \nabla w_h\big)
+ \Big(\big(1-\mbox{$\frac{1}{2}$}(q_I^n+q_P^n)\big)(\mathcal{C}^n_h-\mathcal{C}^n), w_h\Big)  \nonumber \\
&\quad - \mbox{$\frac12$} (({\bf{U}}^{n-1}_h-{\bf{U}}^{n-1})\cdot\nabla \mathcal{C}^n_h, w_h)
            + \mbox{$\frac12$} ((q^{n-1}_I-q^{n-1}_P) (\mathcal{C}^n_h-\mathcal{C}^n), w_h)\nonumber \\
&\quad + \mbox{$\frac12$} (({\bf{U}}^{n-1}_h-{\bf{U}}^{n-1}) \cdot \nabla w_h, \mathcal{C}^n_h )
            + 
            ({\bf{U}}^{n-1}\cdot \nabla w_h, \mathcal{C}^n_h-\mathcal{C}^n),
\qquad\forall\, w_h\in S^1_h ,
\end{align}
with the initial condition $\theta_h^0=0$.
From \eqref{auxeq} and \eqref{fauxeq} we see that $\theta_h^n-\theta^n$ satisfies the equation
\begin{align} \label{fauxeq222}
(\gamma D_{\tau}(\theta_h^{n}-\theta^n), w_h)
+((D({\bf{U}}^{n-1})\nabla(\theta_h^{n}-\theta^n), \nabla w_h)
+(\theta_h^{n}-\theta^n, w_h) =0,\nonumber\\
 \forall\, w_h\in S^1_h .
\end{align}
Similarly, subtracting \eqref{fauxeq} and \eqref{tdis1} from  
\eqref{ReWrFD-C} gives
\begin{align} \label{erreqcdis}
(\gamma D_{\tau}(\mathcal{C}^{n}_h-\theta^{n}_h-\mathcal{C}^{n}), w_h)
+(D({\bf {U}}^{n-1})\nabla(\mathcal{C}^{n}_h-\theta^{n}_h-\mathcal{C}^{n}),  
\nabla w_h)
+(\mathcal{C}^{n}_h-\theta^{n}_h-\mathcal{C}^{n}, w_h) =0  ,  
\nonumber \\
\forall\, w_h\in S^1_h . 
\end{align}
Here $\mathcal{C}^{n}_h-\theta^{n}_h$ and $\theta_h^n$ can be viewed as  
finite element approximations of $\mathcal{C}^{n}$ and $\theta^n$,  
respectively. 
In view of \eqref{regu3},
$D({\bf {U}}^{n-1})$ can be viewed as the value of a piecewise 
linear function (in time) at time $t_{n-1}$ 
and therefore, the conditions \eqref{ellipticity}-\eqref{W1infty-condition}  
are satisfied. 
Applying Lemma \ref{lemma2} to \eqref{erreqcdis}  and \eqref{fauxeq222} yields
\begin{align} \label{c}
&\|D_{\tau}(\mathcal{C}^{n}_h-\theta^{n}_h-{\bf {P}}_h 
\mathcal{C}^{n})\|_{L^p(\widetilde{W}^{-1, q})}
+\|\mathcal{C}^{n}_h-\theta^{n}_h-{\bf{P}}_h\mathcal{C}^{n}\|_{L^p(W^{1, q})}  
\nonumber \\
&\le C(\|\mathcal{C}^{n}-{\bf{R}}_h\mathcal{C}^{n}\|_{L^p(W^{1,q})} 
+h^{-1}\|{\bf P}_h\mathcal{C}^0-\mathcal{C}^0_h\|_{L^q}) \nonumber
\\
&\le Ch\|\mathcal{C}^{n}\|_{L^p(W^{2,q})} +Ch\|\mathcal{C}^0\|_{W^{2,q}} 
,
\quad n=1,\dots,m.  \qquad\mbox{ (use \eqref{ellep}, \eqref{ellep-2}  
and \eqref{Lag-inter}) } 
\end{align}
and 
\begin{align}
&\|D_{\tau}(\theta^{n}_h-\theta^{n})\|_{L^p(\widetilde{W}^{-1, q})}+\|\theta^{n}_h-\theta^{n}\|_{L^p(W^{1, q})}\nonumber\\
&\le C\|D_{\tau}(\theta^n_h-{\bf P}_h\theta^n)\|_{L^p(\widetilde{W}^{-1, q})}+C\|\theta^{n}_h-{\bf P}_h\theta^{n}\|_{L^p(W^{1, q})} \nonumber\\
&\quad 
+C\|D_{\tau} \theta^n-{\bf P}_hD_{\tau} \theta^n \|_{L^p(\widetilde{W}^{-1, q})}+C\|\theta^{n}-{\bf P}_h\theta^{n}\|_{L^p(W^{1, q})} \nonumber \\
&\le 
C\|\theta^{n}-{\bf{R}}_h\theta^{n}\|_{L^p(W^{1, q})}+C\|D_{\tau}\theta^{n}\|_{L^p(\widetilde{W}^{-1,q})}               +C\|\theta^n\|_{L^p(W^{1,q})} \nonumber\\
&\le
C\|D_{\tau}\theta^{n}\|_{L^P(\widetilde{W}^{-1,q})}               +C\|\theta^n\|_{L^p(W^{1,q})}, 
\quad n=1,\dots,m,
\end{align}
where we have used \eqref{ellep-2} to derive the last inequality, and \eqref{ellep} to get the second last inequality (with $m=\ell_0=1$ and the dual case $m=\ell_0=-1$). Therefore, 
\begin{align} \label{thetaerr}
\|D_{\tau} \theta^{n}_h\|_{L^p(\widetilde{W}^{-1, q})}+\|\theta^{n}_h\|_{L^p(W^{1, q})}
 &\le
 C(\|D_{\tau}\theta^{n}\|_{L^p(\widetilde{W}^{-1, q})} + \|\theta^{n}\|_{L^p(W^{1, q})}) .
\end{align}
Applying Lemma \ref{lemma1} to \eqref{auxeq} leads to 
\begin{align} \label{theta1q}
&\|D_{\tau}\theta^{n}\|_{L^p(\widetilde{W}^{-1, q})}  
+ \|\theta^{n}\|_{L^p(W^{1, q})} \nonumber\\
&\le C\|(D({\bf{U}}^{n-1}_h)-D({\bf{U}}^{n-1}))\, 
\nabla \mathcal{C}^{n}_h\|_{L^p(L^{q})}
   + C\|\big(1-\mbox{$\frac{1}{2}$}(q_I^n+q_P^n)\big) 
(\mathcal{C}^n_h-\mathcal{C}^n) \|_{L^p(L^{q})}   \nonumber\\
&\quad +C\|({\bf{U}}^{n-1}_h-{\bf{U}}^{n-1}) 
\cdot\nabla \mathcal{C}^n_h \|_{L^p(L^{q})}
+C\|(q^{n-1}_I-q^{n-1}_P)(\mathcal{C}^n_h-\mathcal{C}^n)\|_{L^p(L^q)}  
\nonumber \\
&\quad +C\|({\bf{U}}^{n-1}_h-{\bf{U}}^{n-1}) \mathcal{C}^n_h  \|_{L^p(L^{q})}
+C\|{\bf{U}}^{n-1}(\mathcal{C}^n_h-\mathcal{C}^n)\|_{L^p(L^{q})} \nonumber \\
&=:I_1^n+I_2^n+I_3^n+I_4^n+I_5^n+I_6^n \, . 
\end{align}
By \eqref{U-Uh_LinftyLinf}-\eqref{U-Uh_LpLq}, we have the estimate 
\begin{align}
I_1^n
&=C\|(D({\bf{U}}^{n-1}_h)-D({\bf{U}}^{n-1}))\, 
\nabla \mathcal{C}^{n}_h\|_{L^p(L^{q})}  \nonumber \\
&\le C \|(D({\bf{U}}^{n-1}_h)-D({\bf{U}}^{n-1}))\, 
\nabla (\mathcal{C}^{n}_h-\mathcal{C}^{n})\|_{L^p(L^{q})}
     +C\|(D({\bf{U}}^{n-1}_h)-D({\bf{U}}^{n-1}))\, 
\nabla \mathcal{C}^{n}\|_{L^p(L^{q})} \nonumber \\
&\le C \|{\bf{U}}^{n-1}_h-{\bf{U}}^{n-1}\|_{L^\infty(L^\infty)}  
\|\nabla (\mathcal{C}^{n}_h-\mathcal{C}^{n})\|_{L^p(L^{q})}
     +C\|{\bf{U}}^{n-1}_h-{\bf{U}}^{n-1}\|_{L^p(L^{q})}  
\|\nabla \mathcal{C}^{n}\|_{L^\infty(L^\infty)} \nonumber \\
&\le Ch^{\frac14} \|\nabla (\mathcal{C}^{n}_h 
-\mathcal{C}^{n})\|_{L^p(L^{q})}
     +C(\|\mathcal{C}^{n-1}_h-\mathcal{C}^{n-1}\|_{L^p(L^{q})}+h)  \, . 
\quad n=1,\dots,m ,
\nn 
\end{align}
Similarly, we get 
\begin{align}
 I_3^n
&=C\|({\bf{U}}^{n-1}_h-{\bf{U}}^{n-1})\,\nabla \mathcal{C}^{n}_h\|_{L^p(L^{q})}  \nonumber \\
&\le Ch^{\frac14} \|\nabla (\mathcal{C}^{n}_h-\mathcal{C}^{n})\|_{L^p(L^{q})} 
       +C(\|\mathcal{C}^{n-1}_h-\mathcal{C}^{n-1}\|_{L^p(L^{q})}+h)  ,  
       \nn \\[5pt]
 I_5^n
&=C\|({\bf{U}}^{n-1}_h-{\bf{U}}^{n-1})\, \mathcal{C}^{n}_h\|_{L^p(L^{q})}  
\le Ch^{\frac14} \|\mathcal{C}^{n}_h-\mathcal{C}^{n}\|_{L^p(L^{q})}
       +C(\|\mathcal{C}^{n-1}_h-\mathcal{C}^{n-1}\|_{L^p(L^{q})}+h)   ,
       \nn 
\end{align}
and also 
\begin{align}
I_2^n + I_4^n + I_6^n 
\le C \|\mathcal{C}^{n}_h-\mathcal{C}^{n}\|_{L^p(L^{q})}  \, . 
\nn 
\end{align}
Substituting the estimates of $I_j^n$, $j=1,\dots,6$, into \eqref{thetaerr}-\eqref{theta1q}, we obtain
\begin{align} \label{thetaerr-2}
&\|D_{\tau} \theta^{n}_h\|_{L^p(\widetilde{W}^{-1, q})}+\|\theta^{n}_h\|_{L^p(W^{1, q})} \nonumber \\
&\le Ch^{\frac14} \|\nabla (\mathcal{C}^{n}_h-\mathcal{C}^{n})\|_{L^p(L^{q})}
     +C\|\mathcal{C}^{n}_h-\mathcal{C}^{n}\|_{L^p(L^{q})}+Ch , \quad 
     n=1,\dots,m, 
\end{align}
which together with \eqref{c} implies
\begin{align}
&\|D_{\tau}(\mathcal{C}^{n}_h- {\bf {P}}_h\mathcal{C}^{n})\|_{L^p(\widetilde{W}^{-1, q})}
+\|\mathcal{C}^{n}_h-{\bf {P}}_h\mathcal{C}^{n}\|_{L^p(W^{1, q})}  \nonumber \\
&\le \|D_{\tau}(\mathcal{C}^{n}_h-\theta^{n}_h-{\bf {P}}_h\mathcal{C}^{n})\|_{L^p(\widetilde{W}^{-1, q})}
+\|\mathcal{C}^{n}_h-\theta^{n}_h-{\bf{P}}_h\mathcal{C}^{n}\|_{L^p(W^{1, q})}  \nonumber \\
&\quad
+ \|D_{\tau} \theta^{n}_h\|_{L^p(\widetilde{W}^{-1, q})}+\|\theta^{n}_h\|_{L^p(W^{1, q})} \nonumber \\
&\le Ch^{\frac14} \|\nabla (\mathcal{C}^{n}_h-\mathcal{C}^{n})\|_{L^p(L^{q})}
     +C\|\mathcal{C}^{n}_h-\mathcal{C}^{n}\|_{L^p(L^{q})}+Ch \nonumber \\
&\le Ch^{\frac14} \|\nabla (\mathcal{C}^{n}_h-{\bf {P}}_h\mathcal{C}^{n})\|_{L^p(L^{q})}
     +C\|\mathcal{C}^{n}_h-{\bf {P}}_h\mathcal{C}^{n}\|_{L^p(L^{q})}+Ch ,
     \,\quad n=1,\dots,m,
\end{align}
where we have used \eqref{ellep} to derive the last inequality. 
When $h\le h_3$ for some $h_3>0$, we can get from above result that 
\begin{align} \label{W1qCnCh}
\|D_{\tau}(\mathcal{C}^{n}_h-{\bf {P}}_h\mathcal{C}^{n})\|_{L^p(\widetilde{W}^{-1, q})}
+\|\mathcal{C}^{n}_h-{\bf {P}}_h\mathcal{C}^{n}\|_{L^p(W^{1, q})}
\le C\|\mathcal{C}^{n}_h-{\bf {P}}_h\mathcal{C}^{n}\|_{L^p(L^{q})}+Ch .
\end{align}
By using \eqref{Lag-inter} and the triangle inequality, we further derive that
\begin{align}
\|D_{\tau}(\mathcal{C}^{n}_h-\Pi_h\mathcal{C}^{n})\|_ 
{L^p(\widetilde{W}^{-1, q})}
+\|\mathcal{C}^{n}_h-\Pi_h\mathcal{C}^{n}\|_{L^p(W^{1, q})}
\le C\|\mathcal{C}^{n}_h-\Pi_h\mathcal{C}^{n}\|_{L^p(L^{q})}+Ch ,\\
     \quad n=1,\dots,m  \, ,   
\nn 
\end{align}
and by Lemma \ref{discr-embed}, 
\begin{align}
\|\mathcal{C}^{n}_h-\Pi_h\mathcal{C}^{n}\|_{L^{\infty}({L}^\infty)}  
&\le C(\|D_{\tau}(\mathcal{C}^{n}_h-\Pi_h\mathcal{C}^{n})\|_ 
{L^p(\widetilde{W}^{-1, q})}
+\|\mathcal{C}^{n}_h-\Pi_h\mathcal{C}^{n}\|_{L^p(W^{1, q})}) \nn   \\
&\le C\|\mathcal{C}^{n}_h-\Pi_h\mathcal{C}^{n}\|_{L^p(L^{q})}+Ch  \nn  \\
&\le \mbox{$\frac{1}{2}$}  
\|\mathcal{C}^{n}_h-\Pi_h\mathcal{C}^{n}\|_{L^\infty(L^{\infty})}  
+ C\|\mathcal{C}^{n}_h-\Pi_h\mathcal{C}^{n}\|_{L^1(L^\infty)} +Ch ,
\,\,\, n=1,\dots,m,
\end{align}
Applying Gronwall's inequality, we see that 
\begin{align}
\|\mathcal{C}^{n}_h-\Pi_h\mathcal{C}^{n}\|_{L^{\infty}({L}^\infty)} \le Ch ,
\qquad n=1,\dots,m.
\end{align}
Finally, using \eqref{ellep}, \eqref{Lag-inter} and the triangle inequality,  
we have
\begin{align}\label{LinftyLinftyCnh}
\begin{aligned}
\|\mathcal{C}^{n}_h-{\bf P}_h\mathcal{C}^{n}\|_{L^{\infty}({L}^\infty)}
&\le \|\mathcal{C}^{n}_h-\Pi_h\mathcal{C}^{n}\|_{L^{\infty}({L}^\infty)}
+\|\Pi_h\mathcal{C}^{n}-{\bf P}_h\mathcal{C}^{n}\|_{L^{\infty}({L}^\infty)}   \\
&\le Ch +Ch\|\mathcal{C}^{n}\|_{L^{\infty}(W^{1,\infty})}  
\le  Ch ,
\quad n=1,\dots,m,
\end{aligned}
\end{align}
which completes  the mathematical induction on \eqref{priest}  
when $h \le h_3$ for some $h_3>0$. 
Consequently, \eqref{LinftyLinftyCnh} holds  
for $m=N$ and \eqref{U-Uh_LinftyLinf} holds for $m=N+1$.

By an inverse inequality and \eqref{LinftyLinftyCnh}, we have
\begin{align}
\|{\bf{P}}_h\mathcal{C}^n-\mathcal{C}^n_h\|_{L^{\infty}(W^{1,\infty})}  
\le Ch^{-1}\|{\bf{P}}_h\mathcal{C}^n-\mathcal{C}^n_h\|_{L^{\infty}(L^\infty)} 
\le C ,
\quad n=1,\dots,N.
\end{align}
and therefore, 
\begin{align*}
&\|{\bf U}_h^n\|_{L^\infty}
\le \|{\bf U}_h^n-{\bf U}^n\|_{L^\infty}+\|{\bf U}^n\|_{L^\infty}
\le Ch^{\frac 14}+C
\le C,
&& n=1,\dots,N,\nonumber\\
&\|\mathcal{C}^n_h\|_{W^{1,\infty}}
\le \|{\bf P}_h\mathcal{C}^n-\mathcal{C}^n_h\|_{W^{1,\infty}}
+\|{\bf P}_h\mathcal{C}^n\|_{W^{1,\infty}}
\le C+\|\mathcal{C}^n\|_{W^{1,\infty}}\le C, && n=1,\dots,N,\nonumber
\end{align*}
where we have used \eqref{U-Uh_LinftyLinf} to estimate  
$\|{\bf U}_h^n-{\bf U}^n\|_{L^\infty}$ 
and \eqref{regu3} for   
$\|{\bf U}^n\|_{L^\infty}$ and $\|\mathcal{C}^n\|_{W^{1,\infty}}$,  
respectively.

The proof of Lemma \ref{Lemma:BD} is completed.  \endproof

\subsection{Proof of (\ref{pqerror})}
Now we turn back to the proof of Theorem 2.1. 
We rewrite the system \eqref{e1}-\eqref{e2} into
\begin{align}
&
\mbox{$
-\nabla\cdot \left(\frac{k(x)}{\mu(c^{n-1})}\nabla p^{n-1}\right)=q_I^{n-1}-q_P^{n-1}, 
$} 
\qquad\qquad\qquad\qquad\qquad\qquad\label{rewrp} \\
&\gamma \partial_{t}c^{n}-\nabla\cdot(D({\bf {u}}^{n-1})\nabla c^{n})
+c^{n}
\mbox{$
=\hat{c}q_I^n
+\left(1-\frac{1}{2}\left(q_I^n+q_P^n\right)\right)c^{n}
$}
\nonumber \\
&\mbox{$
\!\quad\qquad\qquad\qquad\qquad\qquad\qquad\qquad -\frac{1}{2}{\bf{u}}^{n-1}\cdot\nabla c^{n}
-\frac{1}{2}\nabla\cdot ({\bf{u}}^{n-1}\,c^{n})+ E^n,
$}  \label{rewrc}
\end{align}
where
\begin{align}\label{exact-u-discr}
&\mbox{$
{\bf u}^{n-1}=\frac{k(x)}{\mu(c^{n-1})}\nabla p^{n-1} ,
$}
\end{align}
and $E^n$ denotes the truncation error of the linearized scheme, given by
\begin{align*}
E^{n}=&\nabla\cdot((D({\bf u}^{n})-D({\bf u}^{n-1}))\nabla c^n)  +({\bf u}^{n-1}-{\bf u}^n)\cdot\nabla c^n 
-\mbox{$\frac{1}{2}$}((q^n_I-q^n_P)-(q^{n-1}_I-q^{n-1}_P))c^n .
\end{align*}
The regularity assumption \eqref{regu} implies
\begin{align} 
\|E^n\|_{L^p(L^q)}\le C\tau.
\nn 
\end{align}

We subtract \eqref{rewrp} from \eqref{fdis2} to get 
\begin{align} 
&
\mbox{$
\left(\frac{k(x)}{\mu(\mathcal{C}^{n}_h)}\nabla(P^{n}_h 
-\overline{\bf{P}}_hp^{n}),
\nabla v_h\right)
$}
 \nonumber \\
&=
\mbox{$
\left(\frac{k(x)}{\mu(\mathcal{C}^{n}_h)}\nabla(p^{n} 
-\overline{\bf{P}}_h p^{n}),
\nabla v_h\right)+\left(\left(\frac{k(x)}{\mu(c^{n})}
-\frac{k(x)}{\mu(\mathcal{C}^{n}_h)}\right)\nabla p^{n}, \nabla v_h\right)     ,
\quad\forall\, v_h\in \mathring S_h^2 .
$}
\nn 
\end{align}
By Lemma \ref{Lemma:BD} and Lemma \ref{Lemma:W1q-FEM}, 
\begin{align}\label{Phn-Phpn-W1q}
\|P^{n}_h-\overline{\bf{P}}_h p^{n}\|_{W^{1,q}}
&\le 
\mbox{$ 
C\left\|\frac{k(x)}{\mu(\mathcal{C}^{n}_h)}\nabla(p^{n} 
-\overline{\bf{P}}_hp^{n})\right\|_{L^q}
    + C\left\|\left(\frac{k(x)}{\mu(c^{n})}
-\frac{k(x)}{\mu(\mathcal{C}^{n}_h)}\right)\nabla p^{n}\right\|_{L^q} 
$} \nonumber \\
&\le C\|p^{n}-\overline{\bf{P}}_h p^{n}\|_{W^{1,q}}+C\|c^n-\mathcal{C}^{n}_h\|_{L^q} \nonumber \\
&\le Ch^2\|p^{n}\|_{W^{3,q}}+C\|c^n-\mathcal{C}^{n}_h\|_{L^q}  ,
\qquad n=0,1,\dots,N.
\end{align}
Moreover,  subtracting \eqref{exact-u-discr} from \eqref{TDIS4} yields
\begin{align} \label{Un-in-terms-c}
\|{\bf{u}}^{n}-{\bf{U}}^{n}_h\|_{L^q} 
&\le 
\mbox{$ \left\|\frac{k(x)}{\mu(\mathcal{C}^{n}_h)}\nabla(P^{n}_h-p^{n})
+\left(\frac{k(x)}{\mu({\mathcal{C}^{n}_h})}-\frac{k(x)}{\mu(c^{n})}\right)
\nabla p^{n} \right\|_{L^q} 
$} \nn\\
&\le
C\|P^{n}_h-p^n\|_{W^{1,q}} + C\|\mathcal{C}^{n}_h-c^{n}\|_{L^{q}}
\|p^{n}\|_{W^{1, \infty}}
\nn\\
&\le Ch^2\|p^{n}\|_{W^{3,q}}+C\|\mathcal{C}^{n}_h-c^n\|_{L^q} ,
\qquad n=0,1,\dots,N,
\end{align}
where we have used \eqref{Phn-Phpn-W1q} to derive the last inequality.

We take the same approach as used for 
$\|\mathcal{C}^n-\mathcal{C}^{n}_h\|_{L^q} $ in the last subsection to 
estimate $\|c^n-\mathcal{C}^{n}_h\|_{L^q} $. 
We rewrite the finite element system \eqref{fdis1} into
{
\begin{align} \label{ReWrFD-C-F}
& (\gamma D_{\tau}\mathcal{C}_h^{n}, w_h)
+(D({\bf {u}}^{n-1})\nabla \mathcal{C}_h^{n}, \nabla w_h)  
+(\mathcal{C}^n_h, w_h) \\
&=\Big(\hat{c}q_I^n + \big(1-\mbox{$\frac{1}{2}$}(q_I^n+q_P^n)\big) 
c^n, w_h\Big)
-\mbox{$\frac{1}{2}$}({\bf{u}}^{n-1}\cdot\nabla c^n , w_h)  
+ \mbox{$\frac{1}{2}$}({\bf{u}}^{n-1} \cdot \nabla w_h, c^n )
+(E^n, w_h)  \nonumber \\
&\quad+\big((D({\bf {u}}^{n-1})-D({\bf {U}}_h^{n-1})) 
\nabla \mathcal{C}_h^{n}, \nabla w_h\big)
+ \Big(\big(1-\mbox{$\frac{1}{2}$}(q_I^n+q_P^n)\big) 
(\mathcal{C}^n_h-c^n) , w_h\Big)  \nonumber \\
&\quad -\mbox{$\frac{1}{2}$}(({\bf{U}}^{n-1}_h-{\bf{u}}^{n-1}) 
\cdot\nabla \mathcal{C}^n_h, w_h)
   +\mbox{$\frac{1}{2}$}((q^{n-1}_I-q^{n-1}_P) (\mathcal{C}^n_h-c^n), w_h)  
\nn \\
&\quad +\mbox{$\frac{1}{2}$}(({\bf{U}}^{n-1}_h-{\bf{u}}^{n-1}) \cdot \nabla w_h, \mathcal{C}^n_h )
            + ({\bf{u}}^{n-1}\cdot \nabla w_h, \mathcal{C}^n_h-c^n)  - (E^n, w_h),
\qquad\forall\, w_h\in S^1_h . \nonumber
\end{align} 
}
In view of the difference between the right-hand sides of \eqref{rewrc}  
and \eqref{ReWrFD-C-F}, and in order to invoke Lemma 3.2, we define $\chi^n$ to be the solution of an auxiliary parabolic equation:
\begin{align}  \label{auxi1}
&\gamma D_\tau \chi^n -\nabla\cdot(D({\bf{u}}^{n-1})\nabla\chi^n) +\chi^n \nonumber\\
&=
-\nabla\cdot((D({\bf{u}}^{n-1})-D({\bf{U}}^{n-1}_h))
\nabla \mathcal{C}^{n}_h) +\big(1-\mbox{$\frac{1}{2}$}(q_I^n+q_P^n)\big)(\mathcal{C}^n_h-c^n) \nonumber \\
&\quad -\mbox{$\frac{1}{2}$}({\bf{U}}^{n-1}_h-{\bf{u}}^{n-1})\cdot\nabla \mathcal{C}^n_h
             +\mbox{$\frac{1}{2}$}(q^{n-1}_I-q^{n-1}_P) (\mathcal{C}^n_h-c^n)   \nonumber \\
&\quad -\mbox{$\frac{1}{2}$}\nabla\cdot\big(({\bf{U}}^{n-1}_h-{\bf{u}}^{n-1}) \mathcal{C}^n_h\big)
            -\nabla\cdot\big({\bf{u}}^{n-1}(\mathcal{C}^n_h-c^n)\big)  - E^n, 
\end{align}
with the boundary and initial conditions
\begin{align*}
&\mbox{$
-D({\bf{u}}^{n-1})\nabla \chi^n\cdot{\bf {n}}
=-(D({\bf{u}}^{n-1})-D({\bf{U}}^{n-1}_h))
\nabla \mathcal{C}^{n}_h\cdot{\bf {n}}
 -\frac{1}{2}({\bf{U}}^{n-1}_h-{\bf{u}}^{n-1})  \mathcal{C}^n_h  \cdot{\bf {n}}
 $} \nonumber \\
&\,\,\,\quad\qquad\qquad\qquad\qquad -{\bf{u}}^{n-1}(\mathcal{C}^n_h-c^n)  \cdot{\bf {n}} &&\mbox{on}\,\,\,\partial\Omega,\\
&\chi^0=0 &&\mbox{in}\,\,\,\Omega .
\end{align*}
The corresponding finite element approximation of \eqref{auxi1} is defined as: find $\chi^n_h\in S^1_h$, such that
 \begin{align} \label{auxi2}
&(\gamma D_{\tau}\chi^{n}_h, w_h)+(D({\bf{u}}^{n-1})\nabla\chi_h^{n}, \nabla w_h)
+(\chi^{n}_h, w_h) \nonumber  \\
&=\big((D({\bf {u}}^{n-1})-D({\bf {U}}_h^{n-1}))\nabla \mathcal{C}_h^{n}, \nabla w_h\big)
+ \left(\big(1-\mbox{$\frac{1}{2}$}(q_I^n+q_P^n)\big)(\mathcal{C}^n_h-c^n) , w_h\right)  \nonumber \\
&\quad -\mbox{$\frac{1}{2}$}(({\bf{U}}^{n-1}_h-{\bf{u}}^{n-1})\cdot\nabla \mathcal{C}^n_h, w_h)
            +\mbox{$\frac{1}{2}$}((q^{n-1}_I-q^{n-1}_P) (\mathcal{C}^n_h-c^n), w_h) \nonumber \\
&\quad +\mbox{$\frac{1}{2}$}(({\bf{U}}^{n-1}_h-{\bf{u}}^{n-1}) \cdot \nabla w_h, \mathcal{C}^n_h )
            + ({\bf{u}}^{n-1}\cdot \nabla w_h, \mathcal{C}^n_h-c^n) -(E^n, w_h), 
\,\,\,\forall\, w_h\in S^1_h , 
\end{align}
with the initial condition $\chi_h^0=0$. By comparing \eqref{auxi1} and \eqref{auxi2}, we see that
\begin{align} \label{auxi333}
(\gamma D_{\tau}(\chi^{n}_h-\chi^{n}), w_h)+(D({\bf{u}}^{n-1})\nabla(\chi^{n}_h-\chi^{n}), \nabla w_h)
+(\chi^{n}_h-\chi^{n}, w_h) =0,\nonumber \\
 \forall\, w_h\in S^1_h .
\end{align}
Subtracting \eqref{auxi2} and \eqref{rewrc} from \eqref{ReWrFD-C-F} yields
\begin{align} \label{rewrp1}
 (\gamma D_{\tau}(\mathcal{C}_h^{n}-\chi^{n}_h)-\partial_tc^n, w_h)
+(D({\bf {u}}^{n-1})\nabla(\mathcal{C}_h^{n}-\chi^{n}_h-c^n), \nabla w_h) +(\mathcal{C}_h^{n}-\chi^{n}_h-c^n, w_h) = 0 ,
\nonumber \\
\forall\, w_h\in S_h^1.
\end{align}
Again $\mathcal{C}_h^{n}-\chi^{n}_h$ can be viewed as the finite 
element approximation of $c^n$. Then by Lemma \ref{lemma2}, 
\begin{align}  \label{cerr1--0}
&\|\mathcal{C}_h^{n}-\chi^{n}_h-{\bf{P}}_hc^{n}\|_{L^p(L^{q})}  \nonumber      \\
&\le  C \|{\bf{P}}_hc^n-{\bf{R}}_hc^n\|_{L^p(L^q)} +C\|{\bf {P}}_hc_0-\mathcal{C}^0_h\|_{L^q}
+C \|\partial_{tt}c^n\|_{L^p(\widetilde{W}^{-1, q})}\tau\nonumber\\
&\le C\|\mathcal{C}^0_h-c_0\|_{L ^q}  +C(\tau+h^2)
 \qquad \mbox{(use \eqref{ellep}-\eqref{Lag-inter})}    .
\end{align}
Similarly, applying Lemma \ref{lemma2} to \eqref{auxi333} yields
\begin{align}
\|\chi^{n}_h\|_{L^p(L^{q})}
&\le \|\chi_h^n-{\bf{P}}_h\chi^{n}\|_{L^p(L^{q})}
+\|{\bf{P}}_h\chi^{n}\|_{L^p(L^{q})}
\qquad\qquad\qquad\qquad\qquad\, \mbox{(triangle inequality)} \nonumber \\
&\le C(\|{\bf{P}}_h\chi^{n}-{\bf {R}}_h\chi^{n}\|_{L^p(L^{q})} 
+C\|\partial_{tt}\chi\|_{L^p(\widetilde{W}^{-1,q})}\tau) 
+C\|\chi^{n}\|_{L^p(L^{q})} \qquad\,\,\mbox{(use \eqref{lemma-4})}
\nonumber \\
&\le Ch\|\chi^{n}\|_{L^p(W^{1,q})}+C\tau+C\|\chi^{n}\|_{L^p(L^{q})} .
\quad\qquad\qquad\qquad\qquad\qquad\mbox{(use \eqref{ellep}-\eqref{ellep-2})}  
\nonumber
\end{align}
Substituting the last inequality into \eqref{cerr1--0}, we have
\begin{align}
\|\mathcal{C}_h^{n}-{\bf{P}}_hc^{n}\|_{L^p(L^{q})}
&\le C\|\mathcal{C}^0_h-c^0\|_{L^q}  +C(\tau+h^2)
+Ch\|\chi^{n}\|_{L^p(W^{1,q})}+C\tau+C\|\chi^{n}\|_{L^p(L^{q})} \nonumber \\
&\le C\|\mathcal{C}^0_h-c^0\|_{L^q}  +C(\tau+h^2)
+Ch\|\chi^{n}\|_{L^p(W^{1,q})}+C\|\chi^{n}\|_{L^{\infty}(L^{\infty})} ,
\end{align}
and therefore, 
\begin{align}\label{cerr1}
\|\mathcal{C}_h^{n}-c^{n}\|_{L^p(L^{q})}
&\le \|\mathcal{C}_h^{n}-{\bf{P}}_hc^{n}\|_{L^p(L^{q})}
+\|{\bf{P}}_hc^{n}-c^{n}\|_{L^p(L^{q})}  \nonumber \\
&\le C\|\mathcal{C}^0_h-c^0\|_{L^q}  +C(\tau+h^2)
+Ch\|\chi^{n}\|_{L^p(W^{1,q})}+C\|\chi^{n}\|_{L^{\infty}(L^{\infty})} ,
\end{align}
where we have used \eqref{ellep} to estimate 
$\|{\bf{P}}_hc^{n}-c^{n}\|_{L^p(L^{q})} $.

Since $2/p+d/q<1$, there exists $p_0\in(2,p)$ such that $2/p_0+d/q<1$.  
To estimate $\|\chi^{n}\|_{L^{\infty}(L^{\infty})} $, we apply 
Lemma \ref{discr-embed} and Lemma \ref{lemma1} 
to \eqref{auxi1} to get 
\begin{align*}
\|\chi^n\|_{L^{\infty}(L^{\infty})}
&\le C(\|D_\tau \chi^n\|_{L^{p_0}(\widetilde W^{-1,q})} + \|\chi^n\|_{L^{p_0}({W}^{1,q})}) \nonumber \\
&\le
C\|(D({\bf{u}}^{n-1})-D({\bf{U}}^{n-1}_h)) \nabla \mathcal{C}^{n}_h\|_{L^{p_0}(L^q)}
+C\|\big(1-\mbox{$\frac{1}{2}$}(q_I^n+q_P^n)\big)(\mathcal{C}^n_h-c^n) \|_{L^{p_0}(L^q)}\nonumber \\
&\quad +C\|({\bf{U}}^{n-1}_h-{\bf{u}}^{n-1})\cdot\nabla \mathcal{C}^n_h\|_{L^{p_0}(L^q)}
  +C\|(q^{n-1}_I-q^{n-1}_P) (\mathcal{C}^n_h-c^n) \|_{L^{p_0}(L^q)} \nonumber \\
&\quad +C\|({\bf{U}}^{n-1}_h-{\bf{u}}^{n-1}) \mathcal{C}^n_h\|_{L^{p_0}(L^q)}
            +C\|{\bf{u}}^{n-1}(\mathcal{C}^n_h-c^n)\|_{L^{p_0}(L^q)}
            +C\|E^n\|_{L^{p_0}(L^q)} \nonumber \\
&\le
C(\|{\bf U}^{n-1}_h-{\bf u}^{n-1}\|_{L^{p_0}(L^q)}+\|\mathcal{C}^n_h-c^n\|_{L^{p_0}(L^q)}+\|E^{n}\|_{L^{p_0}(L^q)}) \nonumber \\
&\le C\|\mathcal{C}^0_h-c^0\|_{L^q}  +
C\|\mathcal{C}^n_h-c^n\|_{L^{p_0}(L^q)} +C(\tau+h^2),
\end{align*}
where we have used Lemma \ref{Lemma:BD} to estimate 
$\|\nabla\mathcal{C}^{n}_h\|_{L^{\infty}}$ and   
$\|\mathcal{C}^{n}_h\|_{L^{\infty}}$, and \eqref{Un-in-terms-c} 
in deriving the last inequality. Similarly, replacing $p_0$ by $p$ 
in the last inequality yields
\begin{align*}
&(\|D_\tau \chi^n\|_{L^{p}(\widetilde W^{-1,q})} + \|\chi^n\|_{L^{p}({W}^{1,q})})
\le
C\|\mathcal{C}^0_h-c^0\|_{L^q}  +C\|\mathcal{C}^n_h-c^n\|_{L^{p}(L^q)} 
+C(\tau+h^2) .
\end{align*}
By substituting the last two estimates into \eqref{cerr1}, we obtain
\begin{align}
&\|\mathcal{C}^n_h-c^n\|_{L^p(L^q)} \nonumber \\
&\le
C\|\mathcal{C}^0_h-c^0\|_{L^q}  +C(\tau+h^2) 
+Ch\|\mathcal{C}^n_h-c^n\|_{L^{p}(L^q)}
+C\|\mathcal{C}^n_h-c^n\|_{L^{p_0}(L^q)} \nonumber \\
&\le
C\|\mathcal{C}^0_h-c^0\|_{L^q}  +C(\tau+h^2) 
+Ch\|\mathcal{C}^n_h-c^n\|_{L^{p}(L^q)}
+ \frac{1}{2} \|\mathcal{C}^n_h-c^n\|_{L^{p}(L^q)}
+C\|\mathcal{C}^n_h-c^n\|_{L^{1}(L^q)} . \nonumber
\end{align}
When $h \le h_4$ for some $h_4>0$, 
we have 
\begin{align}
&\|\mathcal{C}^n_h-c^n\|_{L^p(L^q)} \le
C\|\mathcal{C}^0_h-c^0\|_{L^q}+ C\|\mathcal{C}^n_h-c^n\|_{L^{1}(L^q)}  +C(\tau+h^2)  .
\end{align}
or equivalently 
\begin{align}
&\bigg(\tau\sum_{n=1}^m\|\mathcal{C}^n_h-c^n\|_{L^q}^p\bigg)^{\frac1p} \le
C\|\mathcal{C}^0_h-c^0\|_{L^q} + C\tau\sum_{n=1}^m\|\mathcal{C}^n_h-c^n\|_{L^q} +C(\tau+h^2)  .
\end{align}
By a similar approach, we can obtain the estimate:
\begin{align} \label{ccpp00}
&\bigg(\tau\sum_{n=k+1}^m\|\mathcal{C}^n_h-c^n\|_{L^q}^p\bigg)^{\frac1p} \le
C\|\mathcal{C}^k_h-c^k\|_{L^q}  + C\tau\sum_{n=k}^m\|\mathcal{C}^n_h-c^n\|_{L^q} +C(\tau+h^2) .
\end{align}
By the generalized Gronwall inequality (Lemma \ref{Lemma:GronW}), 
\begin{align}\label{C-err-esti3}
\|\mathcal{C}^N_h- c^N\|_{L^p(L^q)}
\le C\|\mathcal{C}^0_h- c^0\|_{L^p(L^q)}+C(\tau+h^2)
\le C(\tau+h^2)  .
\end{align}

Finally combining the estimates \eqref{Phn-Phpn-W1q}-\eqref{Un-in-terms-c}  
and \eqref{C-err-esti3}, we obtain the following error estimate when 
$\displaystyle h \le h_{p,q} = \min_{1 \le j \le 4} h_j$ and 
$\displaystyle \tau \le \tau_{p,q} = \min_{1\le j \le 5} \tau_j$,  
\begin{align}
&\!\!\!\!\!\|P^{N}_h-p^{N}\|_{L^p(W^{1,q})}+\|{\bf{U}}^{N}_h-{\bf{u}}^{N}
\|_{L^p(L^q)} +\|\mathcal{C}_h^{N}-c^{N}\|_{L^p(L^q)}
\le C(\tau+h^2)  .
\label{pqerrorff}
\end{align}
Since $q>d$, the inequality above implies \eqref{pqerror}. 
This proves Theorem \ref{MainTHM} in the case $\tau\le \tau_{p,q}$ and $h\le h_{p,q}$.

\subsection{The case $\tau\ge \tau_{p,q}$ or $h\ge h_{p,q}$}
For any $\tau$ and $h$, 
substituting $(v_h,w_h)=(P_h^{n-1},\mathcal{C}_h^n)$ 
into \eqref{fdis2}-\eqref{fdis1} yields
\begin{align}
&\|\nabla P_h^{n-1}\|_{L^2}^2\le \|q_I^{n-1}-q_P^{n-1}\|_{L^2}\|P_h^{n-1}\|_{L^2}
\le \|q_I^{n-1}-q_P^{n-1}\|_{L^2}\|\nabla P_h^{n-1}\|_{L^2}, 
\nn \\
&\mbox{$
D_\tau \left(\frac{\gamma}{2}\|\mathcal{C}_h^n\|_{L^2}^2\right)
\le \frac{\gamma}{4\tau }\|\mathcal{C}_h^n\|_{L^2}^2
+ \frac{\tau}{\gamma}\|\hat c q_I^n\|_{L^\infty}^2 ,
$}
\nn 
\end{align}
which further imply
\begin{align}\label{Uniform-L2PhCh}
\max_{0\le n\le N}
\big(\|P_h^{n}\|_{H^1}\ +\|\mathcal{C}_h^n\|_{L^2}\big)
\le C .
\end{align}

If $\tau\ge \tau_{p,q}$, \eqref{W1qCnCh} still holds for 
$h\le h_{p,q}\le h_3$, which implies that 
\begin{align*} 
&\|D_{\tau}(\mathcal{C}^{n}_h-{\bf {P}}_h\mathcal{C}^{n})\|_{L^p(\widetilde{W}^{-1, q})}
+\|\mathcal{C}^{n}_h-{\bf {P}}_h\mathcal{C}^{n}\|_{L^p(W^{1, q})} \\
&\le C\|\mathcal{C}^{n}_h-{\bf {P}}_h\mathcal{C}^{n}\|_{L^p(L^{q})}+Ch \\
&\le \mbox{$\frac12$}  
\|\mathcal{C}^{n}_h-{\bf {P}}_h\mathcal{C}^{n}\|_{L^p(W^{1,q})} 
+C\|\mathcal{C}^{n}_h-{\bf {P}}_h\mathcal{C}^{n}\|_{L^p(L^{2})}+Ch 
\quad\mbox{(use \eqref{interp-ineq-} here)}\\
&\le \mbox{$\frac12$}\|\mathcal{C}^{n}_h-{\bf {P}}_h\mathcal{C}^{n} 
\|_{L^p(W^{1,q})} 
+C ,
\end{align*}
where the last inequality is due to \eqref{Uniform-L2PhCh}. 
Then we see that 
\begin{align} 
\|\mathcal{C}^{n}_h-c^{n}\|_{L^p(W^{1, q})} 
&\le
\|\mathcal{C}^{n}_h-{\bf P}_h\mathcal{C}^{n}\|_{L^p(W^{1, q})} 
+\|{\bf P}_h\mathcal{C}^{n}-c^{n}\|_{L^p(W^{1, q})}   \nonumber \\
&\le C =C\tau_{p,q}^{-1} \tau_{p,q} \le C\tau_{p,q}^{-1}(\tau+h^2).
\end{align}
On the other hand, \eqref{P-Ph-W1infty} and \eqref{U-Uh_LinftyLinf} imply  
that for $h\le h_{p,q}\le h_2$, 
\begin{align}
&\|{\bf{u}}^{n}-{\bf{U}}^{n}_h\|_{L^\infty} + \|p^{n}-P^{n}_h\|_{W^{1,\infty}} 
\nonumber \\
&\le
\|{\bf{u}}^{n}-{\bf{U}}^{n}\|_{L^\infty} + \|p^{n}-P^{n}\|_{W^{1,\infty}}
+\|{\bf{U}}^{n}-{\bf{U}}^{n}_h\|_{L^\infty} + \|P^{n}-P^{n}_h\|_{W^{1,\infty}} \nonumber \\
&
\le C=C\tau_{p,q}^{-1} \tau_{p,q} \le C\tau_{p,q}^{-1} (\tau+h^2) . 
\end{align}
This proves Theorem \ref{MainTHM} in the case $\tau\ge \tau_{p,q}$ 
and $h\le h_{p,q}$. 

If $h\ge h_{p,q}$, by \eqref{Uniform-L2PhCh} and an inverse inequality, 
we have    
\begin{align}
\max_{0\le n\le N}
\big(\|P_h^{n}\|_{W^{1,q}} +\|\mathcal{C}_h^n\|_{L^q}\big)
&\le 
Ch^{\frac{d}{q}-\frac{d}{2}}\max_{0\le n\le N}
\big(\|P_h^{n}\|_{H^1}\ +\|\mathcal{C}_h^n\|_{L^2}\big) \nonumber \\
&\le 
Ch_{p,q}^{\frac{d}{q}-\frac{d}{2}}\max_{0\le n\le N}
\big(\|P_h^{n}\|_{H^1}\ +\|\mathcal{C}_h^n\|_{L^2}\big) 
\le C,  
\end{align}
and therefore, by noting  
$\|{\bf U}_h^{n}\|_{L^q} \le C\|P_h^{n}\|_{W^{1,q}} \le C$,  
\begin{align}
&\max_{0\le n\le N}
\big(\|p^n-P_h^{n}\|_{W^{1,q}} +\|{\bf u}^n-{\bf U}_h^{n}\|_{L^q}  +\|c^n-\mathcal{C}_h^n\|_{L^q}\big) \nonumber \\
&\le C=Ch_{p,q}^{-2}h_{p,q}^2 \le Ch_{p,q}^{-2} (\tau+h^2)  .
\end{align}
This proves Theorem \ref{MainTHM} in the case $h\ge h_{p,q}$. 
\endproof

\section{Proof of Corollary 2.2} \label{sec:cor}

By using an inverse inequality noting [Lemma \ref{Lemma:BD}, \eqref{cCerr}], 
we can derive from \eqref{fdis1} that 
\begin{align}
\|D_{\tau}\mathcal{C}_h^{n}\|_{L^p(L^q)}
&\le
Ch^{-1}\|D({\bf {U}}_h^{n-1})\nabla \mathcal{C}_h^{n}\|_{L^p(L^q)}
+C\|(q_I^n+q_P^n)\mathcal{C}^n_h\|_{L^p(L^q)} \nonumber \\
&\quad  + C\| {\bf{U}}^{n-1}_h\cdot\nabla \mathcal{C}^n_h\|_{L^p(L^q)}
+Ch^{-1}\|{\bf{U}}^{n-1}_h \mathcal{C}^n_h\|_{L^p(L^q)}
+C\|\hat{c}q_I^n\|_{L^p(L^q)} \nonumber \\
&\le
Ch^{-1}\|{\bf {U}}_h^{n-1}\|_{L^\infty(L^\infty)} 
\|\mathcal{C}_h^{n}\|_{L^\infty(W^{1,\infty})}
+C\|\mathcal{C}_h^{n}\|_{L^\infty(L^{\infty})}
 \nonumber \\
&\quad  + C\|{\bf {U}}_h^{n-1}\|_{L^\infty(L^\infty)} 
\|\mathcal{C}_h^{n}\|_{L^\infty(W^{1,\infty})} 
+Ch^{-1}\|{\bf {U}}_h^{n-1}\|_{L^\infty(L^\infty)} 
\|\mathcal{C}_h^{n}\|_{L^\infty(L^\infty)}
+C  \nonumber \\
&\le Ch^{-1}  
\end{align}
which in turn shows  
$\|D_{\tau}(\mathcal{C}_h^{n}-c^n)\|_{L^p(L^q)}
\le Ch^{-1}$ and 
\begin{align}
\|D_{\tau}(\mathcal{C}_h^{n}-c^n)\|_{L^p(L^q)}
&\le C\tau^{-1} \|\mathcal{C}_h^{n}-c^n\|_{L^p(L^q)}
\le C\tau^{-1}(\tau+h^2)
\le C\tau^{-1} .
\end{align} 
Moreover, by the Sobolev interpolation inequality, we have
\begin{align}
\|\mathcal{C}^N_h- c^N\|_{L^\infty(L^q)}
&\le \|\mathcal{C}^0_h- c^0\|_{L^q}
+C\|\mathcal{C}^N_h- c^N\|_{L^p(L^q)}^{1-\frac1p}
  \|D_\tau(\mathcal{C}^N_h- c^N)\|_{L^p(L^q)}^{\frac1p} \nonumber \\
&\le Ch^2 \|c^0\|_{W^{2,q}} + C(\tau+h^2)^{1-\frac1p}  
\min(\tau^{-1},h^{-1})^{\frac1p} \nonumber \\
&\le Ch^2 + C(\tau^{1-\frac2p} +h^{2-\frac3p}) ,
\end{align} 
where we have used \eqref{pqerror} to estimate  
$\|\mathcal{C}^N_h- c^N\|_{L^p(L^q)}$. Since $p$ can be chosen  
arbitrarily large, combining the above inequality 
and \eqref{Phn-Phpn-W1q}-\eqref{Un-in-terms-c}, 
we obtain \eqref{Linftyerror} immediately and 
the proof of Corollary 2.2 is completed.
\quad \endproof

\section{Numerical results}\label{sec:numer}
\setcounter{equation}{0}
In this section we present numerical results to support our theoretical  
analysis. All the computations are performed by using FreeFEM++ \cite{Fem}.

We consider the equations
\begin{align}
\frac{\partial c}{\partial t}-\nabla\cdot(D({\bf {u}})\nabla c) 
+{\bf{u}}\cdot c=g,\label{ex1}\\
-\nabla\cdot\left(\frac{2}{\mu(c)}\nabla p\right)=f\label{ex2}
\end{align}
in the circular domain $\Omega=\{(x, y): (x-0.5)^2+(y-0.5)^2< 0.5^2\}$,  
with 
$$
{\bf{u}}=-\frac{2}{\mu(c)}\nabla p,\quad \mu(c)=1+c,
\quad
D({\bf {u}})=1+0.1|{\bf{u}}| ,
$$ and an artificially constructed exact solution 
\begin{align}
p=100(x-t)^2e^{-t}, \qquad c=0.5+0.2e^{-t}\cos(x)\sin(y) . 
\end{align}
Substituting this exact solution into the equations \eqref{ex1}-\eqref{ex2} yields the source terms $g$, $f$ and the boundary conditions
\begin{align}\label{bound-initi-condi}
{\bf{u}}\cdot {\bf n}=f_b\quad\mbox{and} \quad D({\bf{u}}) 
\nabla c\cdot {\bf{n}}=g_b
\quad\mbox{on}\,\,\,\partial\Omega .
\end{align}
These are the same type of boundary conditions with given nonzero right-hand sides. 

A quasi-uniform triangulation is made by FreeFEM++ with $M$ nodes 
uniformly distributed on the boundary of the circular domain. 
For simplicity, we denote $h=1/M$. 
We solve the system \eqref{ex1}-\eqref{bound-initi-condi} by the  
proposed method on the quasi-uniform mesh up to time $T=1$.  
The $L^2$ and $L^\infty$  errors of the numerical solutions at time $t=1$ are presented in Table \ref{1} 
with a small fixed time step size $\tau=2^{-14}$ such that the errors  
from time discretization can be negligible in observing the convergence rate 
in the spatial direction. We can see from Table \ref{1} that the proposed method provides  
the accuracy of the optimal order $O(1/M^2)$ for both $\mathcal{C}_h^n$ and 
${\bf U}_h^n$. On the other hand, we present in Table \ref{2} the 
$L^2$ and $L^\infty$ errors of the numerical solutions with a small fixed mesh  
size $h=1/256$ to show the convergence rate in the temporal direction. From 
Table \ref{2}, one can observe clearly that the accuracy of the 
proposed method in time direction is of first order. 
The numerical results are consistent with the analysis given in this paper. 

\begin{table}[h]
\centering
\caption{Errors of numerical solutions in spatial direction 
($\tau=2^{-14}$)}
\label{1}
\begin{tabular}{c|ccccc}
\hline
$ h$
&$\|c^N-\mathcal{C}_h^N\|_{L^{2}}$&$\|{\bf{u}}^N-{\bf{U}}_h^N\|_{L^{2}}$ 
&$\|c^N-\mathcal{C}_h^N\|_{L^{\infty}}$&$\|{\bf{u}}^N-{\bf{U}}_h^N\|_{L^{\infty}}$    \\
 \hline
$1/16$&1.3995E-04   &3.0027E-03  &5.1714E-04   &1.7159E-02   \\

$1/32$&2.8838E-05   &6.9765E-04  &1.4176E-04   &5.2594E-03 \\

$1/64$&7.1872E-06   &1.7068E-04   &3.4551E-05   &1.2412E-03 \\
\hline
 order &2.00 &2.02 &2.03 &2.08 \\
\hline
\end{tabular}
\end{table}\vspace{-10pt}

\begin{table}[h]
  \centering
  \caption{Errors of numerical solutions in time direction 
($h=1/256$)}
\label{2}
\begin{tabular}{c|ccccc}
\hline
$ \tau$&$\|c^N-\mathcal{C}_h^N\|_{L^{2}}$&$\|{\bf{u}}^N-{\bf{U}}_h^N\|_{L^{2}}$  &$\|c^N-\mathcal{C}_h^N\|_{L^{\infty}}$&$\|{\bf{u}}^N-{\bf{U}}_h^N\|_{L^{\infty}}$    \\
 \hline
$1/32$  &4.1618E-04 &6.2041E-04  &2.3635E-03 &2.4287E-03 \\
$1/64$  &1.8478E-04 &2.8533E-04  &1.1310E-03 &1.0462E-03  \\
$\,\,\,1/128$&8.5562E-05 &1.3755E-04  &5.3595E-04 &4.7889E-04  \\
\hline
order &1.06 & 1.06 & 1.07 & 1.12 \\
\hline
\end{tabular}
\end{table}


\section{Conclusion}

In this paper, we have presented an error estimate for the system of PDEs governing miscible displacement in porous media with the Bear--Scheidegger diffusion-dispersion coefficient, which is time-dependent and only ``Lipschitz continuous''. 
The analysis utilizes the discrete maximal $L^p$-regularity of finite element solutions of parabolic equations, which was established in \cite{Li,LS2,LS3} for parabolic equations with Lipschitz continuous coefficients in smooth domains, for time-independent coefficients, time-dependent coefficients with semi-discrete finite element method, and  time-dependent coefficients with fully discrete finite element method, respectively. In these articles (as well as this paper), the domain is assumed to be partitioned into triangles or tetrahedra which fit the boundary $\partial\Omega$ exactly, with possibly curved triangles or tetrahedra near on the boundary.

In the two-dimensional case, the finite element space can be naturally extended (or restricted) to the curved triangle near the boundary. However, in the three-dimensional case, if the boundary faces of the tetrahedra do not exactly lie on $\partial\Omega$ then the curved tetrahedra near the boundary should be specifically constructed instead of being an natural extension of the tetrahedra as in the two-dimensional case. 
For example, for a point $x$ on a boundary face of a tetrahedron one can associate a unique point $y=y(x)\in\partial\Omega$ such that 
\begin{equation*}
y = x + {\bf n}(y) d(x),
\end{equation*}
where ${\bf n}(y)$ is the outward unit normal vector on the point $y\in\partial\Omega$, and $d(x)$ is the signed distance from $x$ to $y$. For $x\in\Omega$ there holds $d(x)>0$, and $x\in\R^d\backslash\Omega$ there holds $d(x)\le 0$. 
Such a transition between the interpolated surface $\partial\Omega_h$ and the exact surface $\partial\Omega$ was introduced as a lift operator in \cite{Dziuk88,DziukElliott_acta}. 
For a tetrahedron $\mathfrak{T}$ with a triangular face $e\subset\partial\Omega_h$, the lift of $e$ onto the smooth boundary $\partial\Omega$ is a curved triangle on $\partial\Omega$. The lift of all such triangles on $\partial\Omega_h$ form a curved triangulation of $\partial\Omega$. 
One can define a region 
\begin{equation*}
\mathfrak{\hat T}=\cup_{x\in e}\{x + \theta\nu(y) d(x): \theta\in[0,1)\} .
\end{equation*}
Then $\mathfrak{\hat T}:=\hat\tau\cup\tau$ is a curved tetrahedron which fit the boundary exactly. 

Such a triangulation with possibly curved tetrahedra on the boundary exists theoretically, as shown above, but is not convenient for practical computation. 
In practical computation, people often replace the original domain $\Omega$ by a triangulated polygonal/polyhedral domain $\Omega_h$. For example, FreeFEM++ solved PDEs in this way. 
Therefore, our numerical example in Section \ref{sec:numer} actually neglects the quadrature error on the boundary triangles (neglecting the quadrature on $\Omega\backslash\Omega_h$). 
This gap between theoretical analysis and practical computation by using FreeFEM++ can possibly be filled in the future by either of the following two approaches:

\begin{enumerate}
[label={(\arabic*)},ref=a\arabic*,topsep=2pt,itemsep=0pt,partopsep=1pt,parsep=1ex,leftmargin=23pt]
\item Instead of assuming that the triangulation fit the boundary exactly, one can use the discrete maximal $L^p$-regularity result established by Kashiwabara and Kemmochi \cite{KK-2018}, who worked on the triangulated domain $\Omega_h$ instead of the original domain $\Omega$. 
In order to apply such results to miscible displacement in porous media, one needs to first extend the result of \cite{KK-2018} to parabolic equations with time-dependent Lipschitz continuous coefficients.  

\item Instead of assuming $\Omega$ to be smooth, one can work on a polygonal/polyhedronal domain directly. However,  the discrete maximal $L^p$-regularity of parabolic equations was only established for the Dirichlet boundary condition so far, see \cite{LiSun2017-MCOM}. 
In order to apply such results to miscible displacement in porous media, one needs to first extend the result of \cite{LiSun2017-MCOM} to the Neumann boundary condition.  
In this case, the error estimates in Theorem \ref{MainTHM} can only be proved for some $q$ depending on the interior angles of the corners and edges, instead of all $q\in(d,\infty)$. 

\end{enumerate}

\newpage

\section*{Appendix: Proof of Lemmas \ref{lemma2}--\ref{Lemma:GronW}}
\renewcommand{\theequation}{A.\arabic{equation}}
$\,$\medskip

{\bf Proof of Lemma 3.2.} \eqref{lemma-3} and \eqref{lemma-4}  can be found in  \cite[(1.18)]{LS3} and \cite[(2.4)]{LS3}, respectively. We prove \eqref{lemma-5} by using \cite[(2.3)]{LS3}, which implies 
(via using inverse inequality) 
\begin{align}\label{Phi-error}
&\|{\bf P}_h\Phi^n-\phi^n_h\|_{L^p(W^{1,q})} \\
&\le
Ch^{-1}\|{\bf P}_h\Phi^n-\phi^n_h\|_{L^p(L^q)}  \nonumber \\
&\le 
Ch^{-1}(\|{\bf P}_h\Phi^n-{\bf R}_h\Phi^n\|_{L^p(L^q)}+\|{\bf P}_h\Phi^0-\phi^0_h\|_{L^q})
&&\mbox{(use \cite[(2.3)]{LS3})} \nonumber\\
&\le 
Ch^{-1}\|\Phi^n-{\bf R}_h\Phi^n\|_{L^p(L^q)}+Ch^{-1}\|{\bf P}_h\Phi^0-\phi^0_h\|_{L^q} 
&&\mbox{(use $L^q$ stability of ${\bf P}_h$)} \nonumber\\
&\le 
C\|\Phi^n-{\bf R}_h\Phi^n\|_{L^p(W^{1,q})}+Ch^{-1}\|{\bf P}_h\Phi^0-\phi^0_h\|_{L^q}.
&&\mbox{(use \eqref{ellep-2} with $l=0$)}\nonumber
\end{align}
From (3.5) and (3.6) we derive 
\begin{align}
&(D_{\tau}({\bf P}_h\Phi^n-\phi^n_h), v_h)+(a(\cdot,t)\nabla({\bf P}_h\Phi^n-\phi^n_h), \nabla v_h)+({\bf P}_h\Phi^n-\phi^n_h, v_h)\\
=&(a(\cdot,t)\nabla({\bf P}_h\Phi^n-{\bf R}_h\Phi^n), \nabla v_h) , 
\quad n=1,\dots,N, \nonumber
\end{align}
which implies
\begin{align}
\|D_{\tau}({\bf P}_h\Phi^n-\phi^n_h)\|_{L^p(\widetilde{W}^{-1,q})} \nonumber
\le&
C\|a(\cdot,t)\nabla({\bf P}_h\Phi^n-\phi^n_h)\|_{L^p(L^q)}+C\|{\bf P}_h\Phi^n-\phi^n_h\|_{L^p(L^q)}\\
&+C\|a(\cdot,t)\nabla({\bf P}_h\Phi^n-{\bf R}_h\Phi^n)\|_{L^p(L^q)}\nonumber\\
\le&C\|{\bf P}_h\Phi^n-\phi^n_h\|_{L^p(W^{1,q})}+C\|{\bf P}_h(\Phi^n-{\bf R}_h\Phi^n)\|_{L^p(W^{1,q})}\nonumber\\
\le&
C\|\Phi^n-{\bf R}_h\Phi^n\|_{L^p(W^{1,q})}+Ch^{-1}\|{\bf P}_h\Phi^0-\phi^0_h\|_{L^q} , 
\end{align}
where we have used \eqref{Phi-error} in the last inequality. 
The proof is completed. \endproof\medskip

{\bf Proof of Lemma \ref{lemma3.3}.}

(1)
Under the conditions of Lemma \ref{lemma3.3}, the Lax--Milgram lemma implies that \eqref{Neumann-Elliptic} has a unique weak solution $u\in H^1\hookrightarrow L^6$ under the constraint $\int_\Omega u\d x=0$, satisfying $\|u\|_{H^1}\le C\|f\|_{L^2}$. Thus $u$ is also a weak solution of 
\begin{align}\label{Neumann-Elliptic-2}
\left\{\begin{aligned}
&\sum_{i,j=1}^d\frac{\partial}{\partial x_i}\bigg(a_{ij}\frac{\partial u}{\partial x_j}\bigg) -u=f-u&&\mbox{in}\,\,\,\Omega,\\
&\sum_{i,j=1}^d a_{ij}n_i\partial_ju=0 &&\mbox{on}\,\,\,\partial\Omega,
\end{aligned}\right. 
\end{align}
which satisfies the following estimate (applying \cite[Theorem 2.4.2.7]{Grisvard} with $p=2$) 
\begin{align}
\|u\|_{H^{2}}\le C\|f-u\|_{L^2} 
&\le C(\|f\|_{L^2} +\|u\|_{L^2}) \nonumber \\
&\le C(\|f\|_{L^2} +\|u\|_{H^1}) \nonumber \\
&\le C \|f\|_{L^2}   .
\end{align}
Since $H^2\hookrightarrow L^\infty$ in both two- and three-dimensional spaces, we have 
\begin{align}
\|u\|_{L^\infty}\le \|u\|_{H^{2}}\le C \|f\|_{L^2} .
\end{align}
Applying \cite[Theorem 2.4.2.7]{Grisvard} again yields  
\begin{align}
\|u\|_{W^{2,q}}
&\le C_q\|f+u\|_{L^q} \nonumber \\
&\le C_q(\|f\|_{L^q} +\|u\|_{L^q}) \nonumber \\
&\le C_q(\|f\|_{L^q} +\|u\|_{L^\infty})  \nonumber \\
&\le C_q(\|f\|_{L^q} +\|f\|_{L^2})  \nonumber \\
&\le C_q\|f\|_{L^q}  .
\end{align}
This proves \eqref{W2q-Elliptic}. 

(2) 
By choosing $q>d$ we have $f\in C^{\alpha}\hookrightarrow L^q$. 
\eqref{W2q-Elliptic} implies $u\in W^{2,q}\hookrightarrow C^{1,\alpha}\hookrightarrow  C^{\alpha}$ with $\alpha=1-d/q\in(0,1)$. Thus $u$ is also a solution of 
\begin{align}\label{Neumann-Elliptic-C2}
\left\{\begin{aligned}
&\sum_{i,j=1}^d\frac{\partial}{\partial x_i}\bigg(a_{ij}\frac{\partial u}{\partial x_j}\bigg) -u=f-u \in C^\alpha &&\mbox{in}\,\,\,\Omega,\\
&\sum_{i,j=1}^d a_{ij}n_i\partial_ju=0 &&\mbox{on}\,\,\,\partial\Omega,
\end{aligned}\right. 
\end{align}
which satisfies the following H\"older estimate (applying \cite[Theorem 4.40 and Corollary 4.41]{Lieberman2013}) 
\begin{align}
\|u\|_{C^{2,\alpha}}\le C\|f-u\|_{C^{\alpha}} 
&\le C (\|f\|_{C^{\alpha}} +\|u\|_{C^{\alpha}}) 
\le C\|f\|_{C^{\alpha}} .
\end{align}
This completes the proof of Lemma \ref{C2alpha-Elliptic}. \endproof\medskip

{\bf Proof of Lemma \ref{C1alpha-Elliptic}}.$\,\,$ 
Since $f\in C^{\alpha}\hookrightarrow L^2$, the Lax--Milgram lemma implies the existence of a unique weak solution $u\in \mathring H^{1}\hookrightarrow L^6$, and the $W^{1,s}$ estimate of elliptic equations (cf. \cite[Theorem 1]{AQ}) implies $\|u\|_{W^{1,d+1}}\le C\|f\|_{L^{d+1}} .$ 
Since $W^{1,d+1}\hookrightarrow L^\infty$, it follows that 
\begin{align}\label{uLinfty-fCalpha}
\|u\|_{L^\infty}\le C(\|g\|_{L^\infty}+\|f\|_{L^{d+1}})\le C(\|g\|_{L^\infty}+\|f\|_{C^{\alpha}} ) .
\end{align}

Let $\chi=\chi(t)$ be a smooth cut-off function defined for $t\in[0,2]$ such that $\chi(t)=1$ for $t\in[1,2]$ and $\chi(0)=0$, satisfying $|\partial_t\chi|\le C$. Then $\chi u$ satisfies the parabolic equation ($u$ is time-independent)
\begin{align}
\left\{\begin{aligned}
&\partial_t (\chi u)-\sum_{i,j=1}^d\frac{\partial}{\partial x_i}\bigg(a_{ij}\frac{\partial (\chi u)}{\partial x_j}\bigg) = 
u\partial_t \chi +\chi g+ \sum_{i=1}^d\partial_i (\chi f_i) &&\mbox{in}\,\,\,\Omega\times[0,2],\\
&\sum_{i,j=1}^d a_{ij}n_i\frac{\partial (\chi u)}{\partial x_j}=\sum_{i=1}^d \chi  f_in_i &&\mbox{on}\,\,\,\partial\Omega\times [0,2],\\
&\chi(0)u(x,0)=0 &&\mbox{for}\,\,\,x\in\Omega .
\end{aligned}\right. 
\end{align}
\cite[Theorem 4.30]{Lieberman} immediately implies 
\begin{align}\label{chi-u-C1alpha}
&\|\chi u\|_{L^\infty(0,2;C^{1,\alpha})} \nonumber \\
&\le C\|u\partial_t\chi\|_{L^\infty(0,2;L^\infty)} 
+C\|\chi g\|_{L^\infty(0,2;L^\infty)} 
+C(\|\chi f\|_{L^\infty(0,2;C^{\alpha})} + \|\chi f\|_{C^{\alpha}(0,2;L^\infty)} ) \nonumber \\
&\le C\|u\|_{L^\infty} +C\|g\|_{L^\infty} 
+C\|f\|_{C^{\alpha}} \nonumber \\
&\le C(\|g\|_{L^\infty} +\|f\|_{C^{\alpha}})  ,
\end{align}
where the last inequality is due to \eqref{uLinfty-fCalpha}. 
Since $\chi$ is independent of the $x$ variable and $u$ is independent of the $t$ variable, it follows that 
$$
\|\chi u\|_{L^\infty(0,2;C^{1,\alpha})}=\|\chi \|_{L^\infty} \|u\|_{C^{1,\alpha}} .
$$
Thus \eqref{chi-u-C1alpha} implies 
\begin{align}
\|u\|_{C^{1,\alpha}} \le C (\|g\|_{L^\infty}+\|f\|_{C^{\alpha}} ) .
\end{align}
This completes the proof of Lemma \ref{C1alpha-Elliptic}. \endproof\medskip

{\bf Proof of Lemma \ref{Lemma:W1q-FEM}}$\,\,$ The existence and uniqueness of solution $u_h\in \mathring S_h^r$ is standard. It suffices to prove the estimate \eqref{W1q-FEM-est}. Note that \eqref{FEM-Shr} is equivalent to 
\begin{align}
\big(a\nabla u_h ,\nabla v_h\big)  + (u_h,v_h)
=({\bf f},\nabla v_h) + (u_h,v_h) ,\quad\forall\, v_h\in \mathring S_h^r .
\end{align}
Let $u\in H^1$ be the solution of the PDE problem 
\begin{align}\label{PDE-m-u}
\left\{\begin{aligned}
&-\nabla\cdot (a\nabla u)+ u  
= -\nabla\cdot {\bf f} + u_h   &&\mbox{in}\,\,\,\Omega ,\\
&a\nabla u\cdot{\bf n}={\bf f}\cdot{\bf n} &&\mbox{on}\,\,\,\partial\Omega,
\end{aligned}\right.
\end{align}
so that $u_h$ is the Ritz projection of $u$. Then the $W^{1,q}$ stability of Ritz projections (as an interpolation  \cite[Corollary A.6]{Gei2}) says that 
\begin{align}
\|u_h\|_{W^{1,q}} \le C\|u\|_{W^{1,q}} ,
\end{align}
and the $W^{1,q}$ estimate of elliptic equations (cf. \cite[Theorem 1]{AQ}) says that 
\begin{align}
\|u\|_{W^{1,q}} 
\le C_q(\|{\bf f}\|_{L^q}+\|u_h\|_{L^q})
\le C_q\|{\bf f}\|_{L^q}+C_{q,\epsilon}\|u_h\|_{L^2} +\epsilon\|u_h\|_{W^{1,q}},
\end{align}
where $\epsilon\in(0,1)$ can be arbitrarily small at the expense of enlarging the constant $C_{q,\epsilon}$. The two estimates above imply 
\begin{align}
\begin{aligned}
\|u_h\|_{W^{1,q}} 
&\le C_q\|{\bf f}\|_{L^q}+C_q\|u_h\|_{L^2} \\
&\le C_q\|{\bf f}\|_{L^q}+C_q\|u_h\|_{H^1} \\
&\le C_q\|{\bf f}\|_{L^q}+C_q\|{\bf f}\|_{L^2}\\
&\le C_q\|{\bf f}\|_{L^q}.
\end{aligned}
\end{align}
\endproof\medskip

{\bf Proof of Lemma \ref{discr-embed}}$\,\,$ 
If we define 
\begin{align}\label{Def-ect}
\phi(t)=\left\{
\begin{aligned}
&\frac{t_k-t}{\tau}\phi^{k-1} + \frac{t-t_{k-1}}{\tau}\phi^{k}  , &&\mbox{for}\,\,\, t\in[t_{k-1},t_k],\,\,\, k=1,\dots,n,\\
&\phi(2t_n-t), &&\mbox{for}\,\,\, t\in[t_n,2t_n], \\
&0, &&\mbox{for}\,\,\, t\in[2t_n,\infty), 
\end{aligned}
\right.
\end{align}
then the function $\phi$ is piecewise linear in time and supported in the time interval $[0,2t_n]$, satisfying the following estimate: 
\begin{align}\label{phi-111}
&\|\partial_t\phi\|_{L^p({\mathbb R}_+;\widetilde{W}^{-1, q})} 
+\|\phi\|_{L^p({\mathbb R}_+;W^{1, q})}  
 \le C(\|D_{\tau}\phi^n\|_{L^p(\widetilde{W}^{-1, q})} 
+\|\phi^{n}\|_{L^p(W^{1, q})})  .
\end{align} 
Let $E$ denote a global extension operator which maps $W^{1,q}$ boundedly into $W^{1,q}({\mathbb R}^d)$ and maps $\widetilde W^{-1,q}$ boundedly into $W^{-1,q}({\mathbb R}^d)$, such that $Eu=u$ in $\Omega$ for all $u\in \widetilde W^{-1,q}$. Such an extension operator exists, by reflecting the function with respect to the boundary $\partial\Omega$; see \cite[Theorems 5.19 and 5.22]{Adams}. By the real interpolation method, we have 
\begin{align}\label{Besov-Linfty}
\begin{aligned}
&\mbox{$E$ maps $(\widetilde{W}^{-1, q},W^{1, q})_{1-1/p,p}$ boundedly into $(W^{-1,q}({\mathbb R}^d),W^{1, q}({\mathbb R}^d))_{1-1/p,p}$} \, ,\\
&(W^{-1,q}({\mathbb R}^d),W^{1, q}({\mathbb R}^d))_{1-1/p,p}=B^{1-2/p,q;p}({\mathbb R}^d)\hookrightarrow C^\alpha({\mathbb R}^d),\,\,\,\mbox{for}\,\,\,\alpha\in(0,1-2/p-d/q),
\end{aligned}
\end{align}
where $B^{1-2/p,q;p}({\mathbb R}^d)$ denotes the Besov space in ${\mathbb R}^d$ (cf. \cite[\textsection 7.32]{Adams}), with the embedding property $B^{1-2/p,q;p}({\mathbb R}^d)\hookrightarrow C^{\alpha}({\mathbb R}^d)$ for $0<\alpha<1-2/p-d/q$ (cf. \cite[\textsection 7.34]{Adams}). 
Then the inhomogeneous Sobolev embedding (see \cite[Proposition 1.2.10]{Lunardi95})  
\begin{align}
\|\phi\|_{L^\infty({\mathbb R}_+; (\widetilde{W}^{-1, q}),W^{1, q}))_{1-1/p,p})}
\le C (\|\partial_t\phi\|_{L^p({\mathbb R}_+;\widetilde{W}^{-1, q})} +\|\phi\|_{L^p({\mathbb R}_+;W^{1, q})})  ,
\end{align}
together with \eqref{phi-111}-\eqref{Besov-Linfty}, implies 
\begin{align}\label{phi222}
\begin{aligned}
\|\phi\|_{L^\infty({\mathbb R}_+; C^\alpha)}
&\le C(\|D_{\tau}\phi^n\|_{L^p(\widetilde{W}^{-1, q})} 
+\|\phi^{n}\|_{L^p(W^{1, q})})   .
\end{aligned}
\end{align}
This proves \eqref{ineq-discr-embed}. 

The inequality \eqref{ineq-discr-embed2} can be proved similarly in view of the interpolation result 
\begin{align}\label{Besov-Linfty2}
\begin{aligned}
&(L^{q}({\mathbb R}^d),W^{2, q}({\mathbb R}^d))_{1-1/p,p}=B^{2-2/p,q;p}({\mathbb R}^d)\hookrightarrow C^{1,\alpha}({\mathbb R}^d),\,\,\,\mbox{for}\,\,\,\alpha\in(0,1-2/p-d/q) . 
\end{aligned}
\end{align}

The proof of Lemma \ref{discr-embed} is complete. \endproof\medskip

{\bf Proof of Lemma \ref{Lemma:GronW}}$\,\,$ H\"{o}lder's inequality implies that 
\begin{align*}
\bigg(\tau\sum_{n=k+1}^m|Y^n|^p\bigg)^{\frac1p} 
&\le \alpha \bigg(Y^{k}+\tau\sum_{n=k+1}^m|Y^n| \bigg)+\beta \nonumber \\
&\le \alpha \bigg(Y^{k}+(t_m-t_k)^{1-\frac{1}{p}}\bigg(\tau\sum_{n=k+1}^m|Y^n|^p\bigg)^{\frac1p} \bigg)+\beta .
\end{align*} 
If $(t_m-t_k)^{1-\frac{1}{p}}\le (2\alpha)^{-1}$ then the last inequality is reduced to 
\begin{align}\label{Gronw-cond-22}
\bigg(\tau\sum_{n=k+1}^m|Y^n|^p\bigg)^{\frac1p} 
&\le 2 \alpha  Y^{k-1} + 2\beta .
\end{align} 

Let $\tau_p= \frac{1}{4(2\alpha)^{1/(1-1/p)}}$ and $m=[\frac{1}{2\tau (2\alpha)^{1/(1-1/p)}}] $ so that $(2m\tau)^{1-\frac1p}\le (2\alpha)^{-1}$, and 
$$
2m\tau
=2\tau \bigg[\frac{1}{2\tau (2\alpha)^{1/(1-1/p)}}\bigg] 
\ge \frac{1}{(2\alpha)^{1/(1-1/p)}} -  2\tau 
\ge \frac{1}{2(2\alpha)^{1/(1-1/p)}} ,\quad\mbox{for}\,\,\,\tau\le\tau_p. 
$$
We choose a sequence $0=t_{n_0}<t_{n_1}<\cdots<t_{n_\ell}=T$ (so $n_\ell=N$) in the following way. 

If $t_{n_j}+2m\tau\ge T$ then we choose $t_{n_{j+1}}=T$. 

If $t_{n_j}+2m\tau<T$ then we choose $t_{n_{j+1}}\in [t_{n_j}+m\tau,t_{n_j}+2m\tau]$ 
such that 
$$
Y^{n_{j+1}}=\min_{{n_j}+m+1\le n\le {n_j}+2m} Y^n. 
$$
Then 
$$
Y^{n_{j+1}}
\le \bigg(\frac{1}{m} \sum_{n={n_j}+m+1}^{n_j+2m}|Y^n|^p\bigg)^{\frac1p} 
= \bigg(\frac{1}{m\tau}\, \tau \sum_{n={n_j}+m}^{n_j+2m}|Y^n|^p\bigg)^{\frac1p} 
\le (m\tau)^{-\frac1p}\,\bigg( \tau \sum_{n=n_j+1}^{n_j+2m}|Y^n|^p\bigg)^{\frac1p} ,
$$
and \eqref{Gronw-cond-22} implies 
\begin{align}\label{Gronw-cond-2}
\bigg( \tau \sum_{n=n_j+1}^{n_j+2m}|Y^n|^p\bigg)^{\frac1p} 
&\le 2 \alpha  Y^{n_j} + 2\beta .
\end{align} 
The last two estimates show that 
\begin{align*}
Y^{n_{j+1}}
&\le 2^{\frac1p}\Delta T^{-\frac1p}\, \bigg( \tau \sum_{n=n_j+1}^{n_j+2m}|Y^n|^p\bigg)^{\frac1p} \\
&\le 2^{1+\frac1p}\Delta T^{-\frac1p}\, \alpha Y^{n_j}+2^{1+\frac1p}\Delta T^{-\frac1p}\, \beta \\
&\le C_{\alpha,p}Y^{n_j}+C_{\alpha,p} \beta . 
\end{align*} 
Iterations of the above two estimates give (the number of iterations is bounded by $2(2\alpha)^{1/(1-1/p)}T$)
\begin{align*}
&\max_{0\le j\le \ell-1} Y^{n_{j}}
\le C_{T,\alpha,p}(Y^{0}+\beta),\\
&\max_{0\le j\le \ell-1}  \bigg(\tau \sum_{n=n_j+1}^{n_j+2m}|Y^n|^p\bigg)^{\frac1p} 
\le C_{T,\alpha,p}(Y^{0}+\beta) ,
\end{align*} 
and applying \eqref{Gronw-cond-2} again yields 
\begin{align*}
&\bigg(\tau \sum_{n=n_{\ell-1}+1}^{n_\ell}|Y^n|^p\bigg)^{\frac1p} 
\le 2 \alpha  Y^{n_{\ell-1}} + 2\beta
\le C_{T,\alpha,p}(Y^{0}+\beta) .
\end{align*} 
Since $\ell\le 1+2(2\alpha)^{1/(1-1/p)}T$ (a bounded number independent of $\tau$), the last two inequalities imply \eqref{Gronw-concl}. This completes the proof of Lemma \ref{Lemma:GronW}. 
\endproof

\end{document}